\RequirePackage{ifpdf}
\ifpdf 
\documentclass[pdftex]{sigma}
\else
\documentclass{sigma}
\fi

\begin{document}
\allowdisplaybreaks
\renewcommand{\PaperNumber}{092}

\FirstPageHeading

\renewcommand{\thefootnote}{$\star$}

\ShortArticleName{Centralizer Construction and Skew Representations}

\ArticleName{A $\boldsymbol{q}$-Analogue of the Centralizer Construction\\
and Skew Representations\\ of the Quantum Af\/f\/ine Algebra\footnote{This paper is a contribution
to the Vadim Kuznetsov Memorial Issue ``Integrable Systems and Related Topics''.
The full collection is available at
\href{http://www.emis.de/journals/SIGMA/kuznetsov.html}{http://www.emis.de/journals/SIGMA/kuznetsov.html}}}

\Author{Mark J. HOPKINS and Alexander I. MOLEV~$^*$}
\AuthorNameForHeading{M.J. Hopkins and A.I. Molev}

\Address{School of Mathematics and Statistics, University of Sydney,
NSW 2006, Australia}

\Email{\href{mailto:markh@maths.usyd.edu.au}{markh@maths.usyd.edu.au},
\href{mailto:alexm@maths.usyd.edu.au}{alexm@maths.usyd.edu.au}}
\URLaddressMarked{\url{http://www.maths.usyd.edu.au/u/alexm}}

\ArticleDates{Received October 14, 2006; Published online December 26, 2006}

\Abstract{We prove an analogue of the Sylvester theorem for the generator
matrices of the quantum af\/f\/ine algebra ${\rm U}_q(\widehat{\mathfrak{gl}}_n)$.
We then use it to give an explicit realization of the
skew representations of the quantum af\/f\/ine algebra.
This allows one to identify them in a~simp\-le way
by calculating their highest weight, Drinfeld
polynomials and the Gelfand--Tsetlin character (or $q$-character).
We also apply the quantum Sylvester theorem to
construct \mbox{a~$q$-analogue} of the Olshanski algebra
as a projective limit of certain centralizers in ${\rm U}_q(\mathfrak{gl}_n)$
and show that this limit algebra contains the $q$-Yangian
as a subalgebra.}

\Keywords{quantum af\/f\/ine algebra;
quantum Sylvester theorem; skew representations}

\Classification{81R10}

\begin{flushright}
\it To the memory of Vadim Kuznetsov
\end{flushright}

\newcommand{\End}{{\rm{End}\ts}}
\newcommand{\Hom}{{\rm{Hom}}}
\newcommand{\ch}{{\rm{ch}\ts}}
\newcommand{\non}{\nonumber}
\newcommand{\wt}{\widetilde}
\newcommand{\wh}{\widehat}
\newcommand{\ot}{\otimes}
\newcommand{\la}{\lambda}
\newcommand{\La}{\Lambda}
\newcommand{\al}{\alpha}
\newcommand{\be}{\beta}
\newcommand{\ga}{\gamma}
\newcommand{\si}{\sigma}
\newcommand{\vp}{\varphi}
\newcommand{\de}{\delta^{}}
\newcommand{\om}{\omega^{}}
\newcommand{\hra}{\hookrightarrow}
\newcommand{\ve}{\varepsilon}
\newcommand{\ts}{\,}
\newcommand{\tpr}{t^{\tss\prime}}
\newcommand{\qin}{q^{-1}}
\newcommand{\tss}{\hspace{1pt}}
\newcommand{\U}{ {\rm U}}
\newcommand{\Y}{ {\rm Y}}
\newcommand{\CC}{\mathbb{C}\tss}
\newcommand{\ZZ}{\mathbb{Z}\tss}
\newcommand{\Z}{\mathbb{Z}}
\newcommand{\A}{\mathcal{A}}
\newcommand{\Pc}{\mathcal{P}}
\newcommand{\Qc}{\mathcal{Q}}
\newcommand{\Tc}{\mathcal{T}}
\newcommand{\Bc}{\mathcal{B}}
\newcommand{\Ec}{\mathcal{E}}
\newcommand{\Ar}{{\rm A}}
\newcommand{\Ir}{{\rm I}}
\newcommand{\Zr}{{\rm Z}}
\newcommand{\gl}{\mathfrak{gl}}
\newcommand{\oa}{\mathfrak{o}}
\newcommand{\spa}{\mathfrak{sp}}
\newcommand{\g}{\mathfrak{g}}
\newcommand{\ka}{\mathfrak{k}}
\newcommand{\p}{\mathfrak{p}}
\newcommand{\sll}{\mathfrak{sl}}
\newcommand{\agot}{\mathfrak{a}}
\newcommand{\qdet}{ {\rm qdet}\ts}
\newcommand{\sdet}{ {\rm sdet}\ts}
\newcommand{\sgn}{ {\rm sgn}\ts}
\newcommand{\Sym}{\mathfrak S}
\newcommand{\fand}{\quad\text{and}\quad}
\newcommand{\Fand}{\qquad\text{and}\qquad}
\newcommand{\antiddots}
{\underset{\displaystyle\cdot\quad\ }
{\overset{\displaystyle\quad\ \cdot}{\cdot}}}
\newcommand{\dddots}
{\underset{\displaystyle\quad\ \cdot}
{\overset{\displaystyle\cdot\quad\ }{\cdot}}}
\newcommand{\xiL}{\xi^{}_{\La}}
\newcommand{\midd}{ {\rm middle} }

\renewcommand{\theequation}{\arabic{section}.\arabic{equation}}

\newtheorem{thm}{Theorem}[section]
\newtheorem{lem}[thm]{Lemma}
\newtheorem{prop}[thm]{Proposition}
\newtheorem{cor}[thm]{Corollary}
\newtheorem{conj}[thm]{Conjecture}

\newtheorem{defin}[thm]{Definition}

\newcommand{\bth}{\begin{thm}}
\renewcommand{\eth}{\end{thm}}
\newcommand{\bpr}{\begin{prop}}
\newcommand{\epr}{\end{prop}}
\newcommand{\ble}{\begin{lem}}
\newcommand{\ele}{\end{lem}}
\newcommand{\bco}{\begin{cor}}
\newcommand{\eco}{\end{cor}}
\newcommand{\bde}{\begin{defin}}
\newcommand{\ede}{\end{defin}}
\newcommand{\bex}{\begin{example}}
\newcommand{\eex}{\end{example}}
\newcommand{\bre}{\begin{remark}}
\newcommand{\ere}{\end{remark}}
\newcommand{\bcj}{\begin{conj}}
\newcommand{\ecj}{\end{conj}}

\newcommand{\bal}{\begin{aligned}}
\newcommand{\eal}{\end{aligned}}
\newcommand{\beq}{\begin{equation}}
\newcommand{\eeq}{\end{equation}}
\newcommand{\ben}{\begin{equation*}}
\newcommand{\een}{\end{equation*}}

\newcommand{\bpf}{\begin{proof}}
\newcommand{\epf}{\end{proof}}

\def\beql#1{\begin{equation}\label{#1}}

\section{Introduction}\label{sec:int}
\setcounter{equation}{0}

The {\it quantum affine algebras\/} form
an important family of quantized enveloping algebras.
They were introduced by Drinfeld and Jimbo in 1985
in the general context of quantum groups
and have since been studied by many authors; see
e.g.\ the books by Chari and Pressley~\cite{cp:gq}
or Klimyk and Schm\"{u}dgen~\cite{ks:qg} for a detailed account of
the origins and applications of the quantum group theory.
A particular $RTT$-presentation of the quantum af\/f\/ine algebras and
their rational counterparts (now known as the {\it Yangians},
due to Drinfeld)
was originally introduced and studied in the work of L.D.~Faddeev
and the Leningrad school of mathematical physicists in the late 70's
and early 80's in relation with the
{\it quantum inverse scattering method\/} (QISM);
see e.g. Reshetikhin, Takhtajan and Faddeev~\cite{rtf:ql}
for the basics of this approach.
The representation theory of the quantum af\/f\/ine algebras and
Yangians, as well as
various applications of the theory,
are discussed in the mentioned above book
\cite[Chapter~12]{cp:gq}. Amongst other applications we note
a surprising connection of the QISM I and QISM II algebras
with the hypergeometric functions
discovered by Kuznetsov and Koornwinder~\cite{kk:gh,k:hf}.

By the results of Olshanski~\cite{o:ea,o:ri},
the $A$ type Yangian is recovered as a certain projective limit
of centralizers in the universal enveloping algebras.
In this paper, we develop a $q$-analogue of the Olshanski
construction where the role of the Yangian is played by
the quantum af\/f\/ine algebra. Let us now discuss this approach in more detail.

Let $n$ and $m$ be integers such that $n>m\geqslant 0$.
Denote by $\Ar_m(n)$ the centralizer
of $\gl_{n-m}$ in the universal enveloping algebra $\U(\gl_n)$,
where $\gl_{n-m}$ is regarded as a natural subalgebra of $\gl_n$.
Due to Olshanski~\cite{o:ea,o:ri}, for any f\/ixed $m$
there exists a chain of natural homomorphisms
\beql{chainA}
\Ar_m(m)\leftarrow\Ar_m(m+1)
\leftarrow\cdots
\leftarrow\Ar_m(n)
\leftarrow\cdots
\eeq
which respect the natural f\/iltrations inherited from the
universal enveloping algebra.
So, one can def\/ine the corresponding
projective limit algebra
\beql{lim}
\Ar_m={\rm lim\ts proj\ts} \Ar_m(n),
\qquad n\to\infty
\eeq
in the category of f\/iltered algebras.
In particular, the algebra $\Ar_0$ is the projective limit
of the centers of $\U(\gl_n)$ and it is
isomorphic to an algebra of polynomials in countably many
variables. Moreover, for $m>0$ one has the isomorphism~\cite{o:ea,o:ri}:
\beql{tensor}
\Ar_m\cong\Ar_0\otimes \Y(\gl_m),
\eeq
where $\Y(\gl_m)$ is the Yangian for the Lie algebra $\gl_m$.
The algebra $\Y(\gl_m)$ is a
deformation of the enveloping algebra
$\U(\gl_m\ot\CC[x])$ in the class of Hopf algebras; see Drinfeld~\cite{d:ha}.
A key part of the proof of the
decomposition \eqref{tensor} is a construction
of algebra homomorphisms
\beql{hom}
\Y(\gl_m)\to\Ar_m(n),\qquad n=m,m+1,\dots
\eeq
compatible with the chain (\ref{chainA}) which def\/ine an embedding
\beql{embedd}
\Y(\gl_m)\hookrightarrow \Ar_m.
\eeq

A similar construction can be applied to the orthogonal and symplectic
Lie algebras which leads to the introduction of the corresponding
twisted Yangians; see Olshanski~\cite{o:ty} and
Molev and Olshanski~\cite{mo:cc}. A super analogue of the centralizer
construction was given by Nazarov~\cite{n:yq}
to obtain the Yangian for the queer Lie superalgebra.
Moreover, the homomorphism \eqref{hom} gives rise to a representation
of the Yangian in the homomorphism space
$\Hom_{\gl_m}\big(L(\mu),L(\la)\big)$ which is called the skew
(or elementary) representation. Such representations
were studied by Cherednik~\cite{c:ni} (mainly in the context of
the quantum af\/f\/ine algebras) and by Nazarov and Tarasov~\cite{nt:ry};
see also \cite{m:yt}.

In this paper we give a $q$-analogue of the centralizer construction
where $\U(\gl_n)$ is replaced by the
quantized enveloping algebra $\U_q(\gl_n)$.
We follow the approach of \cite{m:yt} where an alternative
embedding \eqref{embedd} was constructed. It is based upon a version of the
quantum Sylvester theorem (see Theorem~\ref{thm:qs} below); cf.
Gelfand and Retakh~\cite{gr:dm}, Krob and Leclerc~\cite{kl:mi}.
Here the Yangian $\Y(\gl_m)$
should be replaced with the $q$-Yangian $\Y_q(\gl_m)$ which is
a natural subalgebra of the quantum af\/f\/ine algebra $\U_q(\wh\gl_m)$.
We prove that the $q$-Yangian is embedded into an appropriate
$q$-analogue of the Olshanski algebra $\Ar_m$
(see Theorem~\ref{thm:secoadet} below).
In particular, this gives a~proof of the Poincar\'e--Birkhof\/f--Witt
theorem for the $q$-Yangian, and the same method applies to
provide a new proof of this theorem for the Yangian $\Y(\gl_m)$.
We believe that
a $q$-version of the isomorphism
\eqref{tensor} also takes place, however, we do not have a proof.

Using a homomorphism associated with
the quantum Sylvester theorem, we provide an explicit
construction of the skew
representations of
the quantum af\/f\/ine algebra $\U_q(\wh\gl_n)$.
We generally follow
the approach of \cite{m:yt}, simplifying the arguments
with the use of some observations
of Brundan and Kleshchev~\cite{bk:pp}.
The construction allows us to easily identify
the skew representations by calculating
their highest weights, Drinfeld polynomials
and the Gelfand--Tsetlin characters.
In the Yangian case, the Drinfeld polynomials were calculated
by Nazarov and Tarasov~\cite{nt:ry}; see also
\cite{br:rs,c:ni,nn:pt}.
In a more general context,
the Gelfand--Tsetlin characters were
introduced by Brundan and Kleshchev~\cite{bk:rs}
in their study of representations of the shifted Yangians.
They are analogous to the Yangian characters
of Knight~\cite{k:st} and the $q$-characters
of Frenkel and Reshetikhin~\cite{fr:qc}
for the quantum af\/f\/ine algebras.

\section[Definitions and preliminaries]{Def{}initions and preliminaries}
\label{sec:def}
\setcounter{equation}{0}

We shall use an $R$-matrix presentation of the algebra $\U_q(\gl_n)$
following \cite{j:qu}; see also \cite{ks:qg} for more details.
We f\/ix a complex parameter $q$ which is nonzero and not a root of unity.
Consider the $R$-matrix
\beql{rmatrixc}
R=q\ts\sum_i E_{ii}\ot E_{ii}+\sum_{i\ne j} E_{ii}\ot E_{jj}+
(q-\qin)\sum_{i<j}E_{ij}\ot E_{ji}
\end{equation}
which is an element of $\End\CC^n\ot \End\CC^n$, where
the $E_{ij}$ denote the standard matrix units and the indices run over
the set $\{1,\dots,n\}$. The $R$-matrix satisf\/ies the Yang--Baxter equation
\beql{YBEconst}
R_{12}\ts R_{13}\ts  R_{23} =  R_{23}\ts  R_{13}\ts  R_{12},
\end{equation}
where both sides take values in $\End\CC^n\ot \End\CC^n\ot \End\CC^n$ and
the indices indicate the copies of $\End\CC^n$, e.g.,
$R_{12}=R\ot 1$ etc.

The {\it quantized enveloping algebra\/} $\U_q(\gl_n)$ is generated
by elements $t_{ij}$ and $\bar t_{ij}$ with $1\leqslant i,j\leqslant n$
subject to the relations
\begin{gather}
t_{ij}=\bar t_{ji}=0, \qquad 1 \leqslant i<j\leqslant n,\nonumber\\
t_{ii}\ts \bar t_{ii}=\bar t_{ii}\ts t_{ii}=1,\qquad 1\leqslant i\leqslant n,\label{defrel}\\
R\ts T_1T_2=T_2T_1R,\qquad R\ts \overline T_1\overline T_2=
\overline T_2\overline T_1R,\qquad
R\ts \overline T_1T_2=T_2\overline T_1R.\nonumber
\end{gather}
Here $T$ and $\overline T$ are the matrices
\beql{matrt}
T=\sum_{i,j}t_{ij}\ot E_{ij},\qquad \overline T=\sum_{i,j}
\overline t_{ij}\ot E_{ij},
\eeq
which are regarded as elements of the algebra $\U_q(\gl_n)\ot \End\CC^n$.
Both sides of each of the $R$-matrix relations in \eqref{defrel}
are elements of $\U_q(\gl_n)\ot \End\CC^n\ot \End\CC^n$ and the indices
of~$T$ and~$\overline T$ indicate the copies of $\End\CC^n$ where
$T$ or $\overline T$ acts; e.g. $T_1=T\ot 1$. In terms of the
generators the def\/ining relations between the $t_{ij}$
can be written as
\beql{defrelg}
q^{\de_{ij}}\ts t_{ia}\ts t_{jb}-
q^{\de_{ab}}\ts t_{jb}\ts t_{ia}
=(q-\qin)\ts (\de_{b<a} -\de_{i<j})
\ts t_{ja}\ts t_{ib},
\eeq
where $\de_{i<j}$ equals $1$ if $i<j$ and $0$ otherwise.
The relations between the $\bar t_{ij}$
are obtained by replacing $t_{ij}$ by $\bar t_{ij}$ everywhere in
\eqref{defrelg}. Finally, the relations involving both
$t_{ij}$ and $\bar t_{ij}$ have the form
\beql{defrelg2}
q^{\de_{ij}}\ts \bar t_{ia}\ts t_{jb}-
q^{\de_{ab}}\ts t_{jb}\ts \bar t_{ia}
=(q-\qin)\ts (\de_{b<a}\ts t_{ja}\ts \bar t_{ib} -\de_{i<j}\ts
\ts \bar t_{ja}\ts t_{ib}).
\eeq

It is well known that the algebra $\U_q(\gl_n)$
specializes to $\U(\gl_n)$  as $q\to 1$.
To make this more precise, regard $q$ as a formal variable
and $\U_q(\gl_n)$ as an algebra over $\CC(q)$.
Then set $\A=\CC[q,\qin]$ and consider
the $\A$-subalgebra $\U_{\A}$ of $\U_q(\gl_n)$ generated by
the elements
\beql{taugen1}
\frac{t_{ij}}{q-\qin}\quad\text{for}\quad i> j,
\qquad
\frac{\bar t_{ij}}{q-\qin}\quad\text{for}\quad i< j,
\eeq
and
\beql{taugen2}
\frac{t_{ii}-1}{q-1},\qquad
\frac{\bar t_{ii}-1}{q-1},
\eeq
for $i=1,\dots,n$. Then we have an isomorphism \beql{isom}
\U_{\A}\ot_{\A}\CC\cong \U(\gl_n)
\eeq
with the action of $\A$ on
$\CC$ def\/ined via the evaluation $q=1$; see e.g.\
\cite[Section~9.2]{cp:gq}. Note that the
elements \eqref{taugen1} respectively specialize to the elements $E_{ij}$
and $-E_{ij}$ of $\U(\gl_n)$ while the elements
\eqref{taugen2} specialize to $E_{ii}$
and $-E_{ii}$.

The {\it quantized enveloping algebra\/} $\U_q(\sll_n)$ is def\/ined as
the associative algebra
with generators $k_1,\dots,k_{n-1},k_1^{-1},\dots,k_{n-1}^{-1}$,
$e_1,\dots,e_{n-1}$ and
$f_1,\dots,f_{n-1}$ subject to the def\/ining relations
\begin{gather*}
k_i\tss k_j=k_j\tss k_i, \qquad k^{}_i\tss k_i^{-1}=k_i^{-1}\tss k^{}_i=1, \\
k^{}_i\tss e^{}_j\tss k_i^{-1}=q^{\tss a_{ij}}\tss e^{}_j,\qquad
k^{}_i\tss f^{}_j\tss k_i^{-1}=q^{-a_{ij}} f^{}_j,\\
[e_i,f_j]=\delta_{ij}\ts\frac{k^{}_i-k_i^{-1}}{q-\qin},\\
[e_i,e_j]=[f_i,f_j]=0\qquad\text{if}\ \ |i-j|>1,\\
e^2_i\tss e^{}_{j}-(q+\qin)\tss e^{}_i\tss e^{}_{j}\tss e^{}_i
+e^{}_{j}\tss e^2_i=0\qquad\text{if}\ \ |i-j|=1,\\
f^2_i\tss f^{}_{j}-(q+\qin)\tss f^{}_i\tss f^{}_{j}\tss f^{}_i
+f^{}_{j}\tss f^2_i=0\qquad\text{if}\ \ |i-j|=1,
\end{gather*}
where $[a_{ij}]$ denotes the Cartan matrix associated with the Lie algebra
$\sll_n$ so that its only nonzero entries are $a_{ii}=2$ and
$a_{ij}=-1$ for $|i-j|=1$.

We have an embedding $\U_q(\sll_n)\hra\U_q(\gl_n)$
given by
\begin{gather*}
k_i\mapsto t_{ii}\tss\bar t_{i+1,i+1},\qquad k^{-1}_i\mapsto \bar t_{ii}\tss
t_{i+1,i+1},\\
e_{i}\mapsto -\frac{\bar t_{i,i+1}\ts t_{ii}}{q-\qin},\qquad
f_{i}\mapsto \frac{\bar t_{ii}\ts t_{i+1,i}}{q-\qin}.
\end{gather*}
We shall identify $\U_q(\sll_n)$ with a subalgebra of $\U_q(\gl_n)$
via this embedding.

We let $\U_q^-$, $\U_q^+$ and $\U_q^0$ denote the subalgebras of
$\U_q(\gl_n)$ respectively generated by the $t_{ij}$ with $i>j$,
the $\bar t_{ij}$ with $i<j$ and the $t_{ii},\bar t_{ii}$ with all $i$.
It is implied by \cite[Proposition~9.2.2]{cp:gq} that
multiplication def\/ines an isomorphism of vector spaces
\beql{tenprudec}
\U_q^-\ot\U_q^0\ot\U_q^+\cong \U_q(\gl_n).
\eeq
The products $t_{11}^{m_1}\cdots t_{nn}^{m_n}$ with $m_i\in\ZZ$
form a basis of $\U_q^0$.
Moreover, the monomials
\ben
t^{k_{n,n-1}}_{n,n-1}\ts
\cdots\ts t^{k_{n2}}_{n2}
\cdots\ts t^{k_{32}}_{32}\ts t^{k_{n1}}_{n1}\ts
\cdots\ts t^{k_{21}}_{21}
\een
with non-negative powers $k_{ij}$ form a basis of $\U_q^-$, while
the monomials
\ben
\bar t^{\ts k_{12}}_{12}\cdots\ts \bar t^{\ts k_{1n}}_{1n}\ts
\bar t^{\ts k_{23}}_{23}\cdots\ts \bar t^{\ts k_{2n}}_{2n}\cdots\ts
\bar t^{\ts k_{n-1,n}}_{n-1,n}
\een
with non-negative powers $k_{ij}$ form a basis of $\U_q^+$.

Denote by $\Zr_q$ the center of the algebra
$\U_q(\gl_n)$.
Due to the isomorphism \eqref{tenprudec},
any $z\in\Zr_q$ can be regarded as an element of $\U_q^-\ot\U_q^0\ot\U_q^+$.
Denote by $\chi(z)$ the projection of $z$ to the subalgebra $\U_q^0$
so that $\chi(z)-z$ belongs to the left ideal of $\U_q(\gl_n)$
generated by the elements $\bar t_{ij}$ with $i<j$.
The map
\beql{hchhom}
\chi:\Zr_q\to \U_q^0
\eeq
is an algebra homomorphism called the {\it Harish-Chandra homomorphism\/}.
This homomorphism is injective and its image is the subalgebra
of $\U_q^0$ generated by the symmetric polynomials
in $x_1^2,\dots,x_n^2$ and the polynomial
$x_1^{-1}\cdots\tss x_n^{-1}$,
where $x_i=q^{-i+1}\tss t_{ii}$. Moreover,
those polynomials of total degree zero in $x_1,\dots,x_n$ form
the image of the center
of the subalgebra $\U_q(\sll_n)$ under the Harish-Chandra homomorphism;
see e.g. \cite[Proposition~9.2.5]{cp:gq} and
\cite[Section~6.3.4]{ks:qg}.

Now consider the Lie algebra of Laurent polynomials $\gl_n[\la,\la^{-1}]$
in an indeterminate $\la$. We denote it by
$\wh\gl_n$ for brevity. The {\it quantum affine algebra\/}
$\U_q(\wh\gl_n)$ is a deformation of the universal enveloping algebra
$\U(\wh\gl_n)$.
By def\/inition,
$\U_q(\wh\gl_n)$ has countably many
generators $t_{ij}^{(r)}$ and $\bar t_{ij}^{\ts(r)}$ where
$1\leqslant i,j\leqslant n$ and $r$ runs over nonnegative integers.
They are combined into the matrices
\beql{taff}
T(u)=\sum_{i,j=1}^n t_{ij}(u)\ot E_{ij},\qquad
\overline T(u)=\sum_{i,j=1}^n \bar t_{ij}(u)\ot E_{ij},
\eeq
where $t_{ij}(u)$ and $\bar t_{ij}(u)$ are formal series
in $u^{-1}$ and $u$, respectively:
\beql{expa}
t_{ij}(u)=\sum_{r=0}^{\infty}t_{ij}^{(r)}\ts u^{-r},\qquad
\bar t_{ij}(u)=\sum_{r=0}^{\infty}\bar t_{ij}^{\ts(r)}\ts u^{r}.
\eeq
The def\/ining relations are
\begin{gather}
t_{ij}^{(0)}=\bar t_{ji}^{\ts(0)}=0, \qquad 1 \leqslant i<j\leqslant n,\nonumber\\
t_{ii}^{(0)}\ts \bar t_{ii}^{\ts(0)}=\bar t_{ii}^{\ts(0)}
\ts t_{ii}^{(0)}=1,\qquad 1\leqslant i\leqslant n,\nonumber\\
R(u,v)\ts T_1(u)T_2(v)=T_2(v)T_1(u)R(u,v),\label{defrelaff}\\
R(u,v)\ts \overline T_1(u)\overline T_2(v)=
\overline T_2(v)\overline T_1(u)R(u,v),\nonumber\\
R(u,v)\ts \overline T_1(u)T_2(v)=T_2(v)\overline T_1(u)R(u,v),\nonumber
\end{gather}
where we have used the notation of \eqref{defrel}
and $R(u,v)=R^{\tss q}(u,v)$ is
the {\it trigonometric $R$-matrix\/} given by
\begin{gather}
R^{\tss q}(u,v)= (u-v)\sum_{i\ne j}E_{ii}\ot E_{jj}+(\qin u-q\tss v)
\sum_{i}E_{ii}\ot E_{ii} \nonumber\\
\phantom{R^{\tss q}(u,v)=}+ (\qin-q)\tss u\tss\sum_{i> j}E_{ij}\ot
E_{ji}+ (\qin-q)\tss v\tss\sum_{i< j}E_{ij}\ot E_{ji}.\label{trRm}
\end{gather}
It satisf\/ies the Yang--Baxter equation
\beql{YBE}
R_{12}(u,v)  R_{13}(u,w)R_{23}(v,w) =  R_{23}(v,w) R_{13}(u,w) R_{12}(u,v),
\eeq
where both sides take values in $\End\CC^n\ot \End\CC^n\ot \End\CC^n$ and
the indices indicate the copies of $\End\CC^n$, e.g.,
$R_{12}(u,v)=R(u,v)\ot 1$ etc.; see e.g.~\cite{df:it,fm:ha}
for more details on the structure of $\U_q(\wh\gl_n)$.
Note that it is more common in the literature to use the notation
$\U_q(\wh\gl_n)$ for the centrally extended quantum af\/f\/ine algebra.
In this paper, however, we are mainly concerned with
f\/inite-dimensional representations of the quantum af\/f\/ine algebra
where the central element acts trivially, so that the notation should
not cause a confusion.

The quantized enveloping algebra
$\U_q(\gl_n)$ is a natural subalgebra of $\U_q(\wh\gl_n)$
def\/ined by the embedding
\beql{emb}
t_{ij}\mapsto t_{ij}^{(0)},\qquad \bar t_{ij}\mapsto\bar t_{ij}^{\ts(0)}.
\eeq
Moreover, there is an algebra homomorphism $\U_q(\wh\gl_n)\to \U_q(\gl_n)$
called the {\it evaluation homomorphism\/} def\/ined by
\beql{eval}
T(u)\mapsto T-\overline T\ts u^{-1},\qquad
\overline T(u)\mapsto \overline T-T\ts u.
\eeq

Now we recall the well-known construction of the quantum determinants
for the algebra $\U_q(\wh\gl_n)$; see e.g.~\cite{c:ni,j:qu,rtf:ql}.
Consider
the multiple tensor product
$\U_q(\wh\gl_n)\ot (\End\CC^n)^{\ot\tss r}$ and use the notation of \eqref{defrelaff}.
Then we have the following corollary of \eqref{YBE} and \eqref{defrelaff}:
\beql{fundam}
R(u_1,\dots,u_r)\ts T_1(u_1)\cdots T_r(u_r)=T_r(u_r)
\cdots T_1(u_1)\ts  R(u_1,\dots,u_r),
\eeq
where
\beql{Rlong}
R(u_1,\dots,u_r)=\prod_{i<j}R_{ij}(u_i,u_j),
\eeq
with the product taken in the lexicographical order on the pairs $(i,j)$.
The proof of \eqref{fundam} is exactly the
same as for the Yangians; see e.g.~\cite{mno:yc}.
Furthermore, consider the $q$-permutation operator
$P^{\tss q}\in\End(\CC^n\ot\CC^n)$
def\/ined by
\beql{qperm}
P^{\tss q}=\sum_{i}E_{ii}\ot E_{ii}+ q\tss\sum_{i> j}E_{ij}\ot
E_{ji}+ \qin\sum_{i< j}E_{ij}\ot E_{ji}.
\eeq
For $q=1$ the operator $P^{\tss q}$ turns into the usual
permutation operator
\ben
P=\sum_{i,j}E_{ij}\ot E_{ji}.
\een
The action of symmetric group $\Sym_r$ on the space $(\CC^n)^{\ot\tss r}$
can be def\/ined by setting $s_i\mapsto P^{\tss q}_{s_i}:=
P^{\tss q}_{i,i+1}$ for $i=1,\dots,r-1$,
where $s_i$ denotes the transposition $(i,i+1)$.
If $\si=s_{i_1}\cdots s_{i_l}$ is a~reduced decomposition
of an element $\si\in \Sym_r$ we set
$P^{\tss q}_{\si}=P^{\tss q}_{s_{i_1}}\cdots P^{\tss q}_{s_{i_l}}$.
We note some useful formulas for the action of $P^{\tss q}_{\si}$
on $(\CC^n)^{\ot\tss r}$. Denote by $e_1,\dots,e_n$ the canonical
basis vectors of~$\CC^n$. Then for any indices $a_1<\dots<a_r$
we have
\ben
P^{\tss q}_{\si}(e_{a_1}\ot\cdots\ot e_{a_r})
=q^{l(\si)}\ts e_{a_{\si^{-1}(1)}}\ot\cdots\ot e_{a_{\si^{-1}(r)}},
\een
where $l(\si)$ denotes the length of the permutation $\si$.
This implies a more general formula: for any $\tau\in \Sym_r$ we have
\ben
P^{\tss q}_{\si}(e_{a_{\tau(1)}}\ot\cdots\ot e_{a_{\tau(r)}})
=q^{l(\si\tau^{-1})-l(\tau)}\ts
e_{a_{\tau\si^{-1}(1)}}\ot\cdots\ot e_{a_{\tau\si^{-1}(r)}}.
\een
We denote by $A^q_r$ the $q$-antisymmetrizer
\beql{antisym}
A^q_r=\sum_{\si\in\Sym_r}\sgn\ts\si\cdot P^{\tss q}_{\si}.
\eeq
The above formulas imply that form any $\tau\in \Sym_r$ we have
\ben
A^q_r(e_{a_{\tau(1)}}\ot\cdots\ot e_{a_{\tau(r)}})
=(-q)^{-l(\tau)}\tss A^q_r(e_{a_1}\ot\cdots\ot e_{a_r}).
\een
We have the relation in $\ts\End (\CC^n)^{\ot\tss r}$:
\beql{ranti}
R(1,q^{-2},\dots,q^{-2r+2})=
\prod_{0\leqslant i<j\leqslant r-1}(q^{-2i}-q^{-2j})\ts A^q_r.
\eeq
Now \eqref{fundam} implies
\beql{anitt}
A^q_r\ts T_1(u)\cdots T_r(q^{-2r+2}u)=T_r(q^{-2r+2}u)\cdots T_1(u)\ts  A^q_r
\eeq
which equals
\beql{matelmi}
\sum_{a_i,b_i}{t\ts}^{a_1\cdots\ts a_r}_{b_1\cdots\ts b_r}(u)\ot
E_{a_1b_1}\ot\cdots\ot E_{a_rb_r}
\eeq
for some elements ${t\ts}^{a_1\cdots\ts a_r}_{b_1\cdots\ts b_r}(u)$
of $\U_q(\wh\gl_n)[[u^{-1}]]$
which we call the {\it quantum minors\/}.
They can be given by the following formulas which are immediate from the def\/inition.
If $a_1<\cdots<a_r$ then
\beql{qminorgen}
{t\ts}^{a_1\cdots\ts a_r}_{b_1\cdots\ts b_r}(u)=
\sum_{\si\in \Sym_r} (-q)^{-l(\si)} \cdot t_{a_{\si(1)}b_1}(u)\cdots
t_{a_{\si(r)}b_r}(q^{-2r+2}u),
\eeq
and for any $\tau\in\Sym_r$ we have
\beql{qmsym}
{t\ts}^{a_{\tau(1)}\cdots\ts a_{\tau(r)}}_{b_1\cdots\ts b_r}(u)=
(-q)^{l(\tau)}{t\ts}^{a_1\cdots\ts a_r}_{b_1\cdots\ts b_r}(u).
\eeq
If $b_1<\cdots<b_r$ (and the $a_i$ are arbitrary) then
\beql{qminorgen2}
{t\ts}^{a_1\cdots\ts a_r}_{b_1\cdots\ts b_r}(u)=
\sum_{\si\in \Sym_r} (-q)^{l(\si)} \cdot t_{a_rb_{\si(r)}}(q^{-2r+2}u)\cdots
t_{a_1b_{\si(1)}}(u),
\eeq
and for any $\tau\in\Sym_r$ we have
\beql{qmsym2}
{t\ts}^{a_1\cdots\ts a_r}_{b_{\tau(1)}\cdots\ts b_{\tau(r)}}(u)=
(-q)^{-l(\tau)}{t\ts}^{a_1\cdots\ts a_r}_{b_1\cdots\ts b_r}(u).
\eeq
Moreover, the quantum minor is zero if two top or two bottom indices
are equal.

As the def\/ining relations for the generators $\bar t_{ij}^{\ts(r)}$
have the same matrix form as for the $t_{ij}^{(r)}$,
the above argument can be applied to def\/ine the respective
quantum minors
${\bar t\ts}^{a_1\cdots\ts a_r}_{b_1\cdots\ts b_r}(u)\in\U_q(\wh\gl_n)[[u]]$.
They are given by the same formulas \eqref{qminorgen} and \eqref{qminorgen2},
where the $t_{ij}(u)$ are respectively replaced with $\bar t_{ij}(u)$.
Indeed, we have the relation
\beql{anittbar}
A^q_r\ts \overline T_1(u)\cdots \overline T_r(q^{-2r+2}u)=
\overline T_r(q^{-2r+2}u)\cdots \overline T_1(u)\ts  A^q_r
\eeq
analogous to \eqref{anitt} so that both sides are equal to
\beql{matelmibar}
\sum_{a_i,b_i}{\bar t\ts}^{a_1\cdots\ts a_r}_{b_1\cdots\ts b_r}(u)\ot
E_{a_1b_1}\ot\cdots\ot E_{a_rb_r}.
\eeq
Furthermore, for any indices $i$, $j$ we have the well known
relations which are deduced from~\eqref{fundam}:
\begin{gather}\label{center}
[{t}^{}_{c_id_j}(u),
{t\ts}^{c_1\cdots\ts c_r}_{d_1\cdots\ts d_r}(v)]=0,\qquad
[{t}^{}_{c_id_j}(u),
{\bar t\ts}^{c_1\cdots\ts c_r}_{d_1\cdots\ts d_r}(v)]=0,
\end{gather}
and the same holds with ${t}^{}_{c_id_j}(u)$ replaced by
${\bar t}^{}_{c_id_j}(u)$.

The {\it quantum determinants\/}
of the matrices $T(u)$ and $\overline T(u)$ are respectively
def\/ined by the relations
\beql{qdets}
\qdet T(u)={t\ts}^{1\cdots\ts n}_{1\cdots\ts n}(u),\qquad
\qdet \overline T(u)={\bar t\ts}^{1\cdots\ts n}_{1\cdots\ts n}(u).
\end{equation}
Write
\beql{coeffqdet}
\qdet T(u)=\sum_{k=0}^{\infty} d_k\ts u^{-k},\qquad
\qdet \overline T(u)=\sum_{k=0}^{\infty} \bar d_k\ts u^{k},
\qquad d_k,\bar d_k\in \U_q(\wh\gl_n).
\end{equation}
Due to the property \eqref{center}, the coef\/f\/icients $d_k$ and $\bar d_k$
belong to the center of
the algebra $\U_q(\wh\gl_n)$.
Furthermore, we have the relation $d_0\bar d_0=1$
which is implied by \eqref{defrelaff}.

Note that for any nonzero constant $\al$, the mapping
\beql{scal}
T(u)\mapsto T(\al\tss u),\qquad \overline T(u)\mapsto \overline T(\al\tss u)
\eeq
def\/ines an automorphism of the algebra $\U_q(\wh\gl_n)$.

\section{Quantum Sylvester theorem}\label{sec:qst}
\setcounter{equation}{0}

Suppose that $A=[a_{ij}]$ is a numerical $n\times n$ matrix
and let $1\leqslant m\leqslant n$.
For any indices $i,j=1,\dots,m$ introduce the
minors $c_{ij}$ of $A$ corresponding to the rows $i,m+1,\dots,n$
and columns $j,m+1,\dots,n$ so that
\ben
c_{ij}={a\ts}^{i,m+1,\dots,n}_{j,m+1,\dots,n}.
\een
Let $A^{}_{\Qc\Qc}$ be the
submatrix of $A$ whose rows and columns are numbered
by the elements of the set $\Qc=\{m+1,\dots,n\}$.
The classical Sylvester theorem provides a formula
for the determinant of the matrix $C=[c_{ij}]\tss$:
\ben
\det C=\det A\cdot
\bigl(\tss\det A^{}_{\Qc\Qc}\bigr)\ts^{m-1}.
\een

We shall give a $q$-analogue of this theorem where minors of $A$
are replaced by quantum minors of the matrices $T(u)$ or $\overline T(u)$.
Following \cite{kl:mi}, for its proof
we employ a certain $q$-analogue of the complimentary
minor identity.

We shall also be using the algebra $\U_{\qin}(\wh\gl_n)$.
In order to distinguish its generators from those of the algebra
$\U_q(\wh\gl_n)$
we mark them by the symbol ${}^{\circ}$.
In particular, $t^{\ts\circ}_{ij}(u)$ and $\bar t^{\ts\circ}_{ij}(u)$
will denote
the corresponding generating series
while $T^{\ts\circ}(u)$ and $\overline T^{\ts\circ}(u)$
will stand for the generator matrices.

\bpr\label{prop:qqin}
The mapping
\ben
\om_n: T(u)\mapsto T^{\ts\circ}(u)^{-1},\qquad
\overline T(u)\mapsto \overline T^{\ts\circ}(u)^{-1}
\een
defines an isomorphism $\om_n:\U_q(\wh\gl_n)\to \U_{\qin}(\wh\gl_n)$.
\epr

\bpf
The matrix $T^{\ts\circ}(u)$ satisf\/ies the relation
\ben
R^{\tss\qin}(u,v)\ts T^{\ts\circ}_1(u)T^{\ts\circ}_2(v)=
T^{\ts\circ}_2(v)T^{\ts\circ}_1(u)R^{\tss\qin}(u,v).
\een
Multiplying both sides on the left by the product
${T^{\ts\circ}_1}(u)^{-1}{T^{\ts\circ}_2}(v)^{-1}$
and on the right by the product
${T^{\ts\circ}_2}(v)^{-1}{T^{\ts\circ}_1}(u)^{-1}$, then
conjugating by $P$ and swapping the parameters $u$ and $v$ gives
\ben
PR^{\tss\qin}(v,u)P\tss{T^{\ts\circ}_1}(u)^{-1}{T^{\ts\circ}_2}(v)^{-1}=
{T^{\ts\circ}_2}(v)^{-1}{T^{\ts\circ}_1}(u)^{-1} PR^{\tss\qin}(v,u)P.
\een
Observing that $PR^{\tss\qin}(v,u)P=-R^{\tss q}(u,v)$,
and applying the same argument to the two remaining matrix relations
in \eqref{defrelaff},
we conclude that
$\om_n$ def\/ines a homomorphism. This map is obviously invertible with
the inverse given by
\ben
\omega_n^{-1}:T^{\ts\circ}(u)\mapsto T(u)^{-1},\qquad
\overline T^{\ts\circ}(u)\mapsto \overline T(u)^{-1}
\een
so that $\om_n$ is an isomorphism.
\epf

Suppose now that the integer $m$ is such
that $0\leqslant m\leqslant n$. For any two subsets
\ben
\Pc=\{i_1,\dots,i_m\}
\fand
\Qc=
\{j_1,\dots,j_m\}
\een
of $\{1,\dots,n\}$ with
cardinality $m$ let
\ben
\overline\Pc=
\{i_{m+1},\dots,i_n\}
\fand
\overline\Qc=
\{j_{m+1},\dots,j_n\}
\een
be their set complements in $\{1,\dots,n\}$.
We assume that
\begin{gather*}
i_1<\dots<i_m
\fand
j_1<\dots<j_m\tss,
\qquad
i_{m+1}<\dots<i_n
\fand
j_{m+1}<\dots<j_n\tss.
\end{gather*}
For any $n\times n$ matrix $X$,
we will denote by $X^{}_{\Pc\Qc}$ the submatrix
whose rows and columns are numbered by
the elements of the sets $\Pc$ and $\Qc$ respectively.

Each of the sequences $i_1,\dots,i_n$ and $j_1,\dots,j_n$ above
is a permutation of the sequence $1,\dots,n$. Denote these two permutations
by $\mathbf{i}$ and $\mathbf{j}$ respectively.
If $A$ is a $n\times n$ matrix with complex entries
and $B$ is the inverse matrix, then
the following identity for complementary minors of $A$ and $B$ holds:
\beql{numcomp}
\det A
\cdot
b_{\tss i_{m+1}\dots i_n}^{\tss j_{m+1}\dots j_n}
=
\sgn \mathbf{i}\cdot\sgn \mathbf{j}\cdot
a^{\tss i_1\dots i_m}_{\tss j_1\dots j_m}\tss.
\eeq
Here $a^{\tss i_1\dots i_m}_{\tss j_1\dots j_m}$
is the minor of $A$ corresponding to the submatrix
$A^{}_{\tss\Pc\Qc}\,$, and
$b_{\tss i_{m+1}\dots i_n}^{\tss j_{m+1}\dots j_n}$
is the minor of $B$ corresponding to the submatrix
$B_{\tss\overline\Qc\ts\overline\Pc}\,$.
We now give analogues of the identity \eqref{numcomp}
for the matrices $T(u)$ and $\overline T(u)$.

\bth\label{thm:blocks}
We have the identities
\begin{gather*}
\qdet T(u)\cdot\omega_n^{-1}\Big(\tss
{t^{\ts\circ}}_{\ts i_{m+1}\dots i_n}^{\ts j_{m+1}\dots j_n}(q^{-2n+2}\tss u)\tss\Big)
=
(-q)^{l(\mathbf{j})-l(\mathbf{i})}\cdot
t^{\tss i_1\dots i_m}_{\tss j_1\dots j_m}(u),\\
\qdet \overline T(u)\cdot\omega_n^{-1}\Big(\tss
{{\bar t\ts}^{\ts\circ}}_{\ts i_{m+1}\dots i_n}^{\ts j_{m+1}\dots j_n}(q^{-2n+2}\tss u)\tss\Big)
=
(-q)^{l(\mathbf{j})-l(\mathbf{i})}\cdot
{\bar t\ts}^{\ts i_1\dots i_m}_{\ts j_1\dots j_m}(u),
\end{gather*}
where ${t^{\ts\circ}}_{\ts i_{m+1}\dots i_n}^{\ts j_{m+1}\dots j_n}(v)$
and
${{\bar t\ts}^{\ts\circ}}_{\ts
i_{m+1}\dots i_n}^{\ts j_{m+1}\dots j_n}(v)$
denote the quantum minors in the
quantum affine algebra $\U_{\qin}(\wh\gl_n)$.
\eth

\bpf
Both identities are verif\/ied in the same way so we only consider
the f\/irst one.
By the def\/inition of the quantum determinant,
\beql{qdetan}
\qdet T(u)\,A^q_n=A^q_n\,T_1\cdots T_n,
\eeq
where $T_i=T_i(q^{-2i+2}\tss u)$ for $i=1,\dots,n$. Let us multiply both
sides of \eqref{qdetan} by $T_n^{-1}\cdots T_{m+1}^{-1}$ from the right.
Then \eqref{qdetan}
takes the form
\beql{qdetn}
\qdet T(u)\,A^q_n\,T_n^{-1}\cdots T_{m+1}^{-1}=A^q_n\,T_1\cdots T_m\tss.
\eeq
Now we apply both sides of \eqref{qdetn} to the basis vector
\beql{v=}
e_{j_1}\ot\cdots\ot e_{j_m}
\ot
e_{i_{m+1}}\ot\cdots\ot e_{i_n}
\in (\CC^n)^{\otimes n}.
\eeq
The right hand side gives
\beql{rhsqd}
A^q_n\sum_{a_1,\dots,a_m}
t_{a_1j_{1}}(u)\cdots t_{a_mj_m}(q^{-2m+2}\tss u)\tss
e_{a_1}\ot\cdots\ot e_{a_m}
\ot
e_{i_{m+1}}\ot\cdots\ot e_{i_n}.
\eeq
The summation here can obviously be restricted to the sequences
$a_1,\dots,a_m$ which are permutations of the sequence $i_1,\dots,i_m\tss$.
Consider the vector
$
A^q_n(e_1\ot\cdots \ot e_n).
$
The sum \eqref{rhsqd} is proportional
to this vector, with the coef\/f\/icient
\ben
(-q)^{-l(\mathbf{i})}\cdot
t^{\tss i_1\dots i_m}_{\tss j_1\dots j_m}(u)\tss.
\een
We shall use the notation $\tpr_{ij}(u)$ for the entries of the
matrix $T^{-1}(u)$.
By applying the left hand side of \eqref{qdetn} to the basis vector (\ref{v=}),
we obtain
\ben
A^q_n\sum_{b_{m+1},\dots,b_n}
\tpr_{b_ni_n}(q^{-2n+2}\tss u)\cdots \tpr_{b_{m+1}i_{m+1}}
(q^{-2m}\tss u)\tss
e_{j_{1}}\ot\cdots\ot e_{j_m}
\ot
e_{b_{m+1}}\ot\cdots\ot e_{b_n}.
\een
Here the coef\/f\/icient of
$
A^q_n(e_1\ot\cdots \ot e_n)
$
equals
\ben
(-q)^{-l(\mathbf{j})}\cdot \sum_{\si}(-q)^{-l(\si)}\tss
\tpr_{j_{\si(n)}i_n}(q^{-2n+2}\tss u)\cdots
\tpr_{j_{\si(m+1)}i_{m+1}}(q^{-2m}\tss u),
\een
summed over permutations $\si$ of the set $\{m+1,\dots,n\}$.
This coef\/f\/icient can be written as
\ben
(-q)^{-l(\mathbf{j})}\cdot \omega_n^{-1}\Big(\tss
{t^{\ts\circ}}_{\ i_{m+1}\dots i_n}^{\ j_{m+1}\dots j_n}(q^{-2n+2}\tss u)\tss\Big).
\een
Thus the desired identity follows from \eqref{qdetn}.
\epf

\bco\label{cor:blocks}
We have the identities
\begin{gather*}
\qdet T(u)\cdot\omega_n^{-1}\Big(\qdet T^{\ts\circ}(q^{-2n+2}\tss u)\Big)=1,\\
\qdet \overline T(u)\cdot\omega_n^{-1}
\Big(\qdet \overline T^{\ts\circ}_{}(q^{-2n+2}\tss u)\Big)=1,
\end{gather*}
where $\qdet T^{\ts\circ}(v)$ and $\qdet \overline T^{\ts\circ}(v)$
denote the quantum determinants
in $\U_{\qin}(\wh\gl_n)$.
\eco

For any $0\leqslant m\leqslant n$
introduce the homomorphism $\imath_{m,n}$
\beql{phim}
\imath_{m,n}:\U_{\qin}(\wh\gl_m)\to \U_{\qin}(\wh\gl_n)
\eeq
which takes the coef\/f\/icients
of the series $t^{\circ}_{ij}(u)$ and ${\bar t}^{\ts\circ}_{ij}(u)$
to the respective elements of $\U_{\qin}(\wh\gl_n)$
with the same name.
Consider the composition
\beql{psim}
\psi_{m,n}=\omega_n^{-1}\circ \imath_{m,n} \circ \om_m
\eeq
which is an algebra homomorphism
\ben
\psi_{m,n}: \U_q(\wh\gl_m)\to \U_q(\wh\gl_n).
\een
The action of $\psi_{m,n}$ on the quantum minors
is described in the next lemma.

\ble\label{lem:psiqmin}
We have
\ben
\psi_{m,n}:{t\ts}^{a_1\dots\ts a_r}_{b_1\dots\ts b_r}(u)\mapsto
\big(\ts{t\ts}^{m+1\dots\ts n}_{m+1\dots\ts n}(q^{2n-2m}\tss u)\big)^{-1}\cdot
{t\ts}^{a_1\dots\ts a_r,m+1\dots\ts n}_{b_1\dots\ts b_r,m+1\dots\ts n}
(q^{2n-2m}\tss u),
\een
where $a_i,b_i\in\{1,\dots,m\}$, and the
same formula holds with the minors in the $t_{ij}(u)$ replaced with
the respective minors in the $\bar t_{ij}(u)$.
\ele

\bpf
Due to \eqref{qmsym} and \eqref{qmsym2} we may assume
that $a_1<\dots<a_r$ and $b_1<\dots<b_r$.
By Theorem~\ref{thm:blocks} and Corollary~\ref{cor:blocks} we have
\ben
\om_m:{t\ts}^{a_1\dots\ts a_r}_{b_1\dots\ts b_r}(u)\mapsto
(-q)^{l(\mathbf{a})-l(\mathbf{b})}\cdot
{t^{\ts\circ}}_{\tss 1\dots m}^{\tss 1\dots m}(q^{-2m+2}\tss u)^{-1}\cdot
{t^{\ts\circ}}_{\ts a_{r+1}\dots a_m}^{\ts b_{r+1}\dots b_m}(q^{-2m+2}\tss u),
\een
where $a_{r+1}<\dots<a_m$ and $b_{r+1}<\dots<b_m$ are respective
complementary elements to the sets $\{a_1,\dots,a_r\}$ and
$\{b_1,\dots,b_r\}$ in $\{1,\dots,m\}$, while $\mathbf{a}$ and $\mathbf{b}$
denote the permutations $(a_1,\dots, a_m)$ and $(b_1,\dots, b_m)$
of the sequence $(1,\dots,m)$.
Now, applying the homomorphism $\imath_{m,n}$ we may
regard the quantum minors on the right hand side as series with
coef\/f\/icients in the algebra $\U_{\qin}(\wh\gl_n)$. The proof is completed
by another application of Theorem~\ref{thm:blocks}.
The same argument proves the corresponding
formula for the minors in the $\bar t_{ij}(u)$.
\epf

In particular, we have the following
explicit formulas for the images of the generators of $\U_q(\wh\gl_n)$.

\bco\label{cor:psiacqm}
For any $1\leqslant i,j\leqslant m$, we have
\ben
\psi_{m,n}: t_{ij}(u)\mapsto
\big(\ts{t\ts}^{m+1\dots\ts n}_{m+1\dots\ts n}(q^{2n-2m}\tss u)\big)^{-1}\cdot
{t\ts}^{i,m+1\dots\ts n}_{j,m+1\dots\ts n}
(q^{2n-2m}\tss u),
\een
and the
same formula holds with the image of the $\bar t_{ij}(u)$
where the minors in the $t_{ij}(u)$ are replaced with
the respective minors in the $\bar t_{ij}(u)$.
\eco

The subalgebra of $\U_q(\wh\gl_n)$
generated by the elements $t_{ij}^{(r)}$ and $\bar t_{ii}^{\ts(0)}$
was studied e.g.\ in~\cite{c:ni,nt:yg,rtf:ql}.
We call it the {\it $q$-Yangian\/} and denote by $\Y_q(\gl_n)$.
For any $1\leqslant i,j\leqslant m$, introduce
the following series
with coef\/f\/icients in $\Y_q(\gl_n)$
\ben
t^{\ts\flat}_{ij}(u)={t\ts}^{i,m+1\dots\ts n}_{j,m+1\dots\ts n}(u)
\een
and combine them into the matrix
$T^{\ts\flat}(u)=[\tss t^{\ts\flat}_{ij}(u)]\tss$.
Let $T(u)^{}_{\Qc\Qc}$ be the
submatrix of $T(u)$ whose rows and columns are numbered
by the elements of the set
$\Qc=\{m+1,\dots,n\}$.

The following is an analogue of the Sylvester theorem for
the $q$-Yangian.

\bth\label{thm:qs}
The mapping
\beql{homnatut}
t_{ij}(u)\mapsto t^{\ts\flat}_{ij}(u),\qquad 1\leqslant i,j\leqslant m,
\eeq
defines a homomorphism $\Y_q(\gl_m)\to\Y_q(\gl_n)$. Moreover,
we have the identity
\ben
\qdet T^{\ts\flat}(u)=\qdet T(u)\cdot
\qdet T(q^{-2}\tss u)^{}_{\Qc\Qc}\tss\cdots\tss\qdet
T(q^{-2m+2}\tss u)^{}_{\Qc\Qc}.
\een
\eth

\bpf
By Corollary~\ref{cor:psiacqm}, we have
\beql{psimqmqs}
\psi_{m,n}:t_{ij}(q^{-2n+2m}\tss u)\mapsto
\big(\qdet T(u)^{}_{\Qc\Qc}\big)^{-1}\cdot
t^{\ts\flat}_{ij}(u).
\eeq
Relation~\eqref{center} implies that the coef\/f\/icients
of the series $\qdet T(u)^{}_{\Qc\Qc}$
commute with those of the series
$t^{\ts\flat}_{ij}(v)$.
Hence, since $\psi_{m,n}$ is a homomorphism, we
can conclude by the application of the
automorphism \eqref{scal} that the mapping \eqref{homnatut}
def\/ines a homomorphism.
Furthermore, by Lemma~\ref{lem:psiqmin},
\ben
\psi_{m,n}:{t\ts}^{1\dots \ts m}_{1\dots\ts m}(q^{-2n+2m}\tss u)\mapsto
\big(\qdet T(u)^{}_{\Qc\Qc}\big)^{-1}\cdot\qdet T(u).
\een
On the other hand,
expanding the quantum minor, we obtain from \eqref{psimqmqs} that
\ben
\psi_{m,n}:{t\ts}^{1\dots \ts m}_{1\dots\ts m}(q^{-2n+2m}\tss u)\mapsto
\big(\qdet T(u)^{}_{\Qc\Qc}\cdots\qdet T(q^{-2m+2}\tss u)^{}_{\Qc\Qc}\big)^{-1}
\cdot \qdet T^{\ts\flat}(u),
\een
completing the proof of the desired identity for $\qdet T^{\ts\flat}(u)$.
\epf

We shall also need a dual version of the quantum Sylvester theorem.
Let $m\geqslant 0$ and $n\geqslant 1$.
Instead of the homomorphism
$\imath_{m,n}$ def\/ined in \eqref{phim}, consider another homomorphism
\beql{phimj}
\jmath_m:\U_{\qin}(\wh\gl_n)\to \U_{\qin}(\wh\gl_{m+n})
\eeq
which takes the coef\/f\/icients of the series
$t^{\ts\circ}_{ij}(u)$ and $\bar t^{\ts\circ}_{ij}(u)$ to
the respective coef\/f\/icients of the series
$t^{\ts\circ}_{m+i,m+j}(u)$ and $\bar t^{\ts\circ}_{m+i,m+j}(u)$.
Consider the composition
\beql{psimj}
\phi_m=\omega_{m+n}^{-1}\circ \jmath_m \circ \om_n.
\eeq
Then $\phi_m$ is an algebra homomorphism
\ben
\phi_m: \U_q(\wh\gl_n)\to \U_q(\wh\gl_{m+n}).
\een
The corresponding version
of Lemma~\ref{lem:psiqmin} now takes the form

\ble\label{lem:psiqminj}
We have
\ben
\phi_m:{t\ts}^{a_1\dots\ts a_r}_{b_1\dots\ts b_r}(u)\mapsto
\big(\ts{t\ts}^{1\dots\ts m}_{1\dots\ts m}(q^{2m}\tss u)\big)^{-1}\cdot
{t\ts}^{1\dots\ts m,\ts m+a_1\dots\ts m+a_r}_{1\dots\ts m,\ts m+b_1\dots\ts m+b_r}
(q^{2m}\tss u),
\een
where $a_i,b_i\in\{1,\dots,n\}$,
and the
same formula holds with the minors in the $t_{ij}(u)$ replaced with
the respective minors in the $\bar t_{ij}(u)$.
\ele

For any indices $1\leqslant i,j\leqslant n$ introduce
the following series
with coef\/f\/icients in the $q$-Yangian $\Y_q(\gl_{m+n})$,
\ben
t^{\ts\sharp}_{ij}(u)={t\ts}^{1\dots \ts m,m+i}_{1\dots\ts m,m+j}(u),
\een
and combine them into the matrix
$T^{\ts\sharp}(u)=[\tss t^{\ts\sharp}_{ij}(u)]\tss$.
Let $T(u)^{}_{\Pc\Pc}$ be the
submatrix of $T(u)$ whose rows and columns are numbered
by the elements of the set
$\Pc=\{1,\dots,m\}$. The following theorem is proved in the same way as
Theorem~\ref{thm:qs}.

\bth\label{thm:qssh}
The mapping
\beql{homnatutsh}
t_{ij}(u)\mapsto t^{\ts\sharp}_{ij}(u),\qquad 1\leqslant i,j\leqslant n,
\eeq
defines a homomorphism $\Y_q(\gl_{n})\to\Y_q(\gl_{m+n})$. Moreover,
we have the identity
\ben
\qdet T^{\ts\sharp}(u)=\qdet T(u)\cdot
\qdet T(q^{-2}\tss u)^{}_{\Pc\Pc}\tss\cdots\tss\qdet
T(q^{-2n+2}\tss u)^{}_{\Pc\Pc}.
\een
\eth

\bre\label{rem:anabar}
The obvious analogues of
Theorems~\ref{thm:qs} and \ref{thm:qssh}
for the subalgebras of $\U_q(\wh\gl_n)$ generated by
the elements $\bar t_{ij}^{\ts(r)}$ and $t_{ii}$
can be proved by the same argument.
\ere

We also point out a corollary to be used in Section~\ref{sec:sr}.
It follows from Lemma~\ref{lem:psiqminj} with the use of
the automorphism \eqref{scal}.

\bco\label{cor:imqmin}
The mapping
\begin{gather}
t_{ij}(u) \mapsto
\big(\ts{t\ts}^{1\dots\ts m}_{1\dots\ts m}(u)\big)^{-1}\cdot
{t\ts}^{1\dots\ts m,\ts m+i}_{1\dots\ts m,\ts m+j}(u),\nonumber\\
\bar t_{ij}(u) \mapsto
\big(\ts{\bar t\ts}^{1\dots\ts m}_{1\dots\ts m}(u)\big)^{-1}\cdot
{\bar t\ts}^{1\dots\ts m,\ts m+i}_{1\dots\ts m,\ts m+j}(u),\label{homnatshar}
\end{gather}
defines a homomorphism $\U_q(\wh\gl_{n})\to\U_q(\wh\gl_{m+n})$. Moreover,
the images of the quantum minors
under the homomorphism \eqref{homnatshar} are found by
\ben
{t\ts}^{a_1\dots\ts a_r}_{b_1\dots\ts b_r}(u)\mapsto
\big(\ts{t\ts}^{1\dots\ts m}_{1\dots\ts m}(u)\big)^{-1}\cdot
{t\ts}^{1\dots\ts m,\ts m+a_1\dots\ts m+a_r}_{1\dots\ts m,\ts m+b_1\dots\ts m+b_r}
(u),
\een
where $a_i,b_i\in\{1,\dots,n\}$, and the
same formula holds with the minors in the $t_{ij}(u)$ replaced with
the respective minors in the $\bar t_{ij}(u)$.
\eco

\section[Skew representations of $\U_q(\wh\gl_n)$]{Skew representations of $\boldsymbol{\U_q(\wh\gl_n)}$}
\label{sec:sr}
\setcounter{equation}{0}

Consider the quantum af\/f\/ine algebra $\U_q(\wh\gl_n)$
and suppose, as before, that the complex parameter $q$
is nonzero and not a root of unity.
The irreducible ({\it pseudo}) {\it highest weight representation\/}
$L(\nu(u),\bar\nu(u))$ of $\U_q(\wh\gl_n)$
with the ({\it pseudo}) {\it highest weight\/} $(\nu(u),\bar\nu(u))$
is generated by a nonzero vector $\zeta$ such that
\begin{alignat}{3}
t_{ij}(u)\ts\zeta&=0,\qquad &&\bar t_{ij}(u)\ts\zeta=0 \qquad &&\text{for} \quad
1\leqslant i<j\leqslant n,
\non\\
t_{ii}(u)\ts\zeta&=\nu_i(u)\ts\zeta,\qquad
&&\bar t_{ii}(u)\ts\zeta=\bar \nu_i(u)\ts\zeta
\qquad &&\text{for} \quad 1\leqslant i\leqslant n,
\non
\end{alignat}
where $\nu(u)=(\nu_1(u),\dots,\nu_n(u))$ and
$\bar\nu(u)=(\bar\nu_1(u),\dots,\bar\nu_n(u))$ are certain $n$-tuples
of formal power series in $u^{-1}$ and $u$, respectively.

Suppose that there exist
polynomials $P_1(u),\dots,P_{n-1}(u)$ in $u$,
all with constant term $1$, such that
\beql{nunupone}
\frac{\nu_k(u)}{\nu_{k+1}(u)}=q^{-\deg P_k}\cdot\frac{P_k(u\tss q^2)}{P_k(u)}
=\frac{\bar\nu_k(u)}{\bar\nu_{k+1}(u)}
\eeq
for any $k=1,\dots,n-1$. The f\/irst equality in \eqref{nunupone}
is understood in the sense
that the ratio of polynomials has to be
expanded as a power series in $u^{-1}$,
while for the second equality the same ratio has to be expanded
as a power series in $u$. Then the representation
$L(\nu(u),\bar\nu(u))$ is f\/inite-dimensional.
The $P_k(u)$ are called the {\it Drinfeld polynomials\/} of this
representation. The families of Drinfeld polynomials
parameterize
the type $1$ f\/inite-dimensional irreducible representations
of the subalgebra $\U_q(\wh\sll_n)$ of $\U_q(\wh\gl_n)$;
see \cite[Section~12]{cp:gq} and
\cite{df:it,fm:ha}.

Let $\la=(\la_1,\dots,\la_n)$
be an $n$-tuple of integers such that $\la_1\geqslant\dots\geqslant\la_n$.
The corresponding irreducible highest weight representation
$L(\la)$ of $\U_q(\gl_n)$
is generated by a nonzero vector $\xi$ such that
\begin{gather*}
\bar t_{ij}\ts\xi=0 \quad \text{for} \quad
1\leqslant i<j\leqslant n,
\\
t_{ii}\ts\xi=q^{\la_i}\tss\xi \quad
\text{for} \quad 1\leqslant i\leqslant n.
\end{gather*}
This representation is a $q$-analogue of the irreducible
$\gl_n$-module
with the highest weight $\la$. In particular, these modules have the same
dimension.

We shall need an analogue of the Gelfand--Tsetlin basis
for the module $L(\la)$ \cite{j:qr}; see also
\cite[Section~7.3.3]{ks:qg} for more details.
A {\it pattern\/} $\La$ (associated with $\la$) is a sequence
of rows of integers $\La_n,\La_{n-1},\dots,\La_1$, where
$\La_k=(\lambda_{k1},\dots,\lambda_{kk})$ is the $k$-th row from the bottom,
the top row $\La_n$
coincides with $\la$, and the following
{\it betweenness conditions\/} are satisf\/ied: for $k=1,\dots,n-1$
\beql{betw}
\la_{k+1,i+1}\leqslant\la_{ki}\leqslant \la_{k+1,i}
\qquad{\rm for}\quad i=1,\dots,k.
\end{equation}
We shall be using the notation $l_{ki}=\la_{ki}-i+1$.
Also, for any integer $m$ we set
\beql{int}
[m]=\frac{q^m-q^{-m}}{q-\qin}.
\end{equation}
There exists a basis $\{\xiL\}$ in $L(\la)$ parameterized
by the patterns $\La$ such that the action of the generators of
$\U_q(\gl_n)$ is given by
\begin{gather}\label{tgt}
t_{kk}\ts\xiL=q^{w_k}\ts\xiL,\qquad w_k=\sum_{i=1}^k\la_{ki}-\sum_{i=1}^{k-1}\la_{k-1,i},
\\
\label{egt}
e_k\ts\xiL=-\sum_{j=1}^k\frac{[\ts l_{k+1,1}-l_{kj}]\cdots
[\ts l_{k+1,k+1}-l_{kj}]}
{[\ts l_{k1}-l_{kj}]\cdots\wedge_j\cdots [\ts l_{kk}-l_{kj}]}
\ts \xi^{}_{\La+\delta_{kj}},\\
\label{fgt}
f_k\ts\xiL=\sum_{j=1}^k\frac{[\ts l_{k-1,1}-l_{kj}]\cdots [\ts l_{k-1,k-1}-l_{kj}]}
{[\ts l_{k1}-l_{kj}]\cdots\wedge_j\cdots [\ts l_{kk}-l_{kj}]}
\ts \xi^{}_{\La-\delta_{kj}},
\end{gather}
where $\La\pm\delta_{kj}$ is obtained from
$\La$ by replacing the entry $\lambda_{kj}$ with $\lambda_{kj}\pm1$,
and $\xi^{}_{\La}$ is supposed to be equal to zero
if $\La$ is not a pattern; the symbol $\wedge_j$ indicates
that the $j$-th factor is skipped.

The representation $L(\la)$ of $\U_q(\gl_n)$
can be extended to
a representation of $\U_q(\wh\gl_n)$ via the
evaluation homomorphism; see \eqref{eval}. This $\U_q(\wh\gl_n)$-module
is called the {\it evaluation module\/}. Its Drinfeld
polynomials are given by
\beql{evalmod}
P_k(u)=(1-q^{2\la_{k+1}}u)(1-q^{2\la_{k+1}+2}u)\cdots(1-q^{2\la_{k}-2}u),
\eeq
for $k=1,\dots,n-1$. We shall be using some formulas for the action
of certain quantum minors of the matrices $T(u)$ and $\overline T(u)$
in the Gelfand--Tsetlin basis of $L(\la)$.
Set
\beql{ttil}
T_{ij}(u)=\frac{u\ts t_{ij}-u^{-1}\ts \bar t_{ij}}{q-\qin}.
\eeq
Clearly, $(q-\qin)\ts T_{ij}(u)$ equals $u$ times the image of
$t_{ij}(u^2)$ under the evaluation homomorphism~\eqref{eval}.
We also def\/ine the corresponding
quantum minors
${T\ts}^{a_1\cdots\ts a_r}_{b_1\cdots\ts b_r}(u)$ by using the formulas
\eqref{qminorgen} or \eqref{qminorgen2}, so that, in particular,
if $a_1<\cdots<a_r$ then
\beql{qminortb}
{T\ts}^{a_1\cdots\ts a_r}_{b_1\cdots\ts b_r}(u)=
\sum_{\sigma\in \Sym_r} (-q)^{-l(\sigma)} \cdot T_{a_{\sigma(1)}b_1}(u)\cdots
T_{a_{\sigma(r)}b_r}(q^{-r+1}u).
\eeq
For any $k=1,\dots,n$ we have in $L(\la)$,
\beql{qmactgtd}
{T\ts}^{1\cdots\ts k}_{1\cdots\ts k}(u)\ts\xiL
=\prod_{i=1}^k \frac{u\tss q^{\tss l_{ki}}-u^{-1}q^{-l_{ki}}}{q-\qin}\ts\xiL.
\eeq
Indeed, the coef\/f\/icients of the polynomial
${T\ts}^{1\cdots\ts k}_{1\cdots\ts k}(u)$ commute with the elements
of the subalgebra $\U_q(\gl_k)$ and so the relation
is verif\/ied by the application of the polynomial to the highest vector
of the $\U_q(\gl_k)$-module $L(\La_k)$.
Furthermore, for any $k=1,\dots,n-1$ and $j=1,\dots,k$ we have
\begin{align}\label{qmactgtr}
{T\ts}^{1\cdots\ts k}_{1\cdots\ts k-1,\ts k+1}(q^{-l_{kj}})\ts\xiL
&=-[\ts l_{k+1,1}-l_{kj}]\cdots [\ts l_{k+1,k+1}-l_{kj}]\ts
\xi^{}_{\La+\delta_{kj}}\\
\intertext{and}
\label{qmactgtl}
{T\ts}^{1\cdots\ts k-1,\ts k+1}_{1\cdots\ts k}(q^{-l_{kj}})\ts\xiL
&=[\ts l_{k-1,1}-l_{kj}]\cdots [\ts l_{k-1,k-1}-l_{kj}]\ts
\xi^{}_{\La-\delta_{kj}}.
\end{align}
The formulas \eqref{qmactgtr} and
\eqref{qmactgtl} appeared for the f\/irst time in \cite{nt:yg}
in a slightly dif\/ferent form. For the proof of \eqref{qmactgtr},
it suf\/f\/ices to use the relation
\ben
{T\ts}^{1\cdots\ts k}_{1\cdots\ts k-1,\ts k+1}(u)
={T\ts}^{1\cdots\ts k}_{1\cdots\ts k}(u)\ts e_k-
q\ts e_k\ts {T\ts}^{1\cdots\ts k}_{1\cdots\ts k}(u),
\een
which can be deduced from \eqref{fundam},
and then apply \eqref{egt} and \eqref{qmactgtd}. The proof
of \eqref{qmactgtl} is similar with the use of
\eqref{fgt} and \eqref{qmactgtd}.

Suppose now that $m$ is a non-negative integer and let
$\la=(\la_1,\dots,\la_{m+n})$
be an $(m+n)$-tuple of integers such that
$\la_1\geqslant\dots\geqslant\la_{m+n}$.
Furthermore, let $\mu=(\mu_1,\dots,\mu_m)$
be an $m$-tuple of integers such that
$\mu_1\geqslant\dots\geqslant\mu_m$. Regarding
$\U_q(\gl_m)$ as a natural subalgebra of $\U_q(\gl_{m+n})$,
consider the space $\Hom_{\U_q(\gl_m)}\big(L(\mu),L(\la)\big)$.
This vector space is isomorphic to the subspace
$L(\la)^+_{\mu}$ of $L(\la)$ which consists of $\U_q(\gl_m)$-singular
vectors of weight $\mu$,
\begin{gather*}
L(\la)^+_{\mu}=\{\eta\in L(\la)\ |\ \bar t_{ij}\tss\eta=0
\quad 1\leqslant i<j\leqslant m\Fand\\
\phantom{L(\la)^+_{\mu}=\{\eta\in L(\la)\ |\ }{} t_{ii}\tss \eta=q^{\mu_i}\tss\eta\quad i=1,\dots,m\}.
\end{gather*}
Note that $L(\la)^+_{\mu}$ is nonzero if and only if
\ben
\la_i\geqslant \mu_i\geqslant \la_{i+n}\qquad\text{for}\quad i=1,\dots,m;
\een
see e.g. \cite{cp:gq}. We shall assume that these inequalities hold.
In that case, $L(\la)^+_{\mu}$ admits a basis~$\{\zeta^{}_{\La}\}$ labelled by
the trapezium Gelfand--Tsetlin patterns $\Lambda$ of the form{\samepage
\begin{align}
\lambda^{}_{1}\qquad\qquad\lambda^{}_{2}\qquad\qquad\lambda^{}_{3}
\qquad\qquad\qquad&\cdots\qquad\qquad
\lambda^{}_{m+n-1}\qquad\qquad\lambda^{}_{m+n}\non\\
\qquad\lambda^{}_{m+n-1,1}\qquad\lambda^{}_{m+n-1,2}
\qquad\qquad\qquad&\cdots\qquad\qquad\qquad\lambda^{}_{m+n-1,m+n-1}\non\\
\qquad\qquad\qquad\dddots\qquad\dddots
\qquad\qquad\qquad\quad&\cdots\qquad\qquad\qquad\antiddots\non\\
\qquad\qquad\qquad\qquad\lambda^{}_{m+1,1}\quad \lambda^{}_{m+1,2}
\qquad\quad&\cdots\quad\qquad\lambda^{}_{m+1,m+1}\non\\
\quad\qquad\qquad\qquad\qquad\qquad\mu^{}_{1}\qquad\ \mu^{}_{2}
\qquad&\cdots\qquad\mu^{}_{m}\non
\end{align}}
These arrays are formed by
integers $\la_{ki}$
satisfying the betweenness conditions
\ben
\la^{}_{k+1,i+1}\leqslant\la^{}_{ki}\leqslant\la^{}_{k+1,i}
\een
for $k=m,m+1,\dots, m+n-1$ and $1\leqslant i\leqslant k$,
where we have set
\ben
\la^{}_{mi}=\mu_i,\quad i=1,\dots,m\Fand
\la^{}_{m+n,j}=\la_j,\quad j=1,\dots,m+n.
\een

Observe that $L(\la)^+_{\mu}$ is a natural module over
the subalgebra $\U'_q(\gl_n)$ of $\U_q(\gl_{m+n})$ generated
by the elements $t_{ij}$ and $\bar t_{ij}$ with
$m+1\leqslant i,j\leqslant m+n$. Consider
the trapezium pattern $\La^0$
with the entries given by
\beql{lao}
\la^0_{m+k,i}=\min\{\la_i,\mu_{i-k}\},
\qquad k=1,\dots,n-1,\quad i=1,\dots,m+k,
\eeq
where we assume $\mu_j=+\infty$ if $j\leqslant 0$. Then
\ben
t_{pp}\ts\zeta_{\La^0}=q^{w_p}\ts\zeta_{\La^0},\qquad w_p=
\sum_{i=1}^p\la^0_{pi}-\sum_{i=1}^{p-1}\la^0_{p-1,i},
\qquad p=m+1,\dots,m+n.
\een
The weight $(w_{m+1},\dots,w_{m+n})$ of $\zeta_{\La^0}$
is maximal with respect to the standard ordering on the set of weights
of the $\U'_q(\gl_n)$-module $L(\la)^+_{\mu}$.

We now make $L(\la)^+_{\mu}$ into a module over the quantum af\/f\/ine
algebra $\U_q(\wh\gl_n)$ by taking the composition
of the homomorphism \eqref{homnatshar} with the evaluation
homomorphism
$\U_q(\wh\gl_{m+n})\to \U_q(\gl_{m+n})$
given by \eqref{eval}. Thus,
we can write the formulas for the action of $\U_q(\wh\gl_n)$
on $L(\la)^+_{\mu}$ in the form
\ben
\bal
t_{ij}(u)&\mapsto
\Big(\qdet(T-\overline T\tss u^{-1})^{}_{\A\A}\Big)^{-1}\cdot
\qdet(T-\overline T\tss u^{-1})^{}_{\A_i\A_j},\\
\bar t_{ij}(u)&\mapsto
\Big(\qdet(\overline T-T\tss u)^{}_{\A\A}\Big)^{-1}\cdot
\qdet(\overline T-T\tss u)^{}_{\A_i\A_j},
\eal
\een
where $\A=\{1,\dots,m\}$ and $\A_i=\{1,\dots,m,m+i\}$
while $T$ and $\overline T$ denote the generator matrices
for the algebra $\U_q(\wh\gl_{m+n})$. The subspace $L(\la)^+_{\mu}$
of $L(\la)$ is preserved by these operators
due to relations \eqref{center}.

\bth\label{thm:irred}
The representation $L(\la)^+_{\mu}$ of $\U_q(\wh\gl_n)$
is irreducible.
\eth

\bpf
By Corollary~\ref{cor:imqmin}, for the action of quantum minors
on $L(\la)^+_{\mu}$ we can write
\ben
{t\ts}^{a_1\cdots\ts a_r}_{b_1\cdots\ts b_r}(u^2)
\mapsto
{T\ts}^{1\cdots\ts m}_{1\cdots\ts m}(u)^{-1}\ts
{T\ts}^{1\cdots\ts m,m+a_1\cdots\ts m+a_r}_{1\cdots\ts m,m+b_1\cdots\ts m+b_r}(u)
\prod_{i=1}^r\frac{q-\qin}{u\ts q^{-m-i+1}},
\een
where $a_i,b_i\in\{1,\dots,n\}$. Note that the coef\/f\/icients
of the polynomial ${T\ts}^{1\cdots\ts m}_{1\cdots\ts m}(u)$ act
on $L(\la)^+_{\mu}$ as scalar operators found from
\ben
{T\ts}^{1\cdots\ts m}_{1\cdots\ts m}(u)\ts\zeta^{}_{\La}
=\prod_{i=1}^m
\frac{u\tss q^{\tss \ell_{i}}-u^{-1}q^{-\ell_{i}}}
{q-\qin}\ts\zeta^{}_{\La},\qquad \ell_i=\mu_i-i+1.
\een
On the other hand, by the formulas \eqref{qmactgtr} and \eqref{qmactgtl},
each of the operators
\ben
{T\ts}^{1\cdots\ts m,m+1\cdots\ts m+k}_{1\cdots\ts m,
m+1\cdots\ts m+k-1,m+k+1}(q^{-l_{m+k,i}})\fand
{T\ts}^{1\cdots\ts m,m+1\cdots\ts m+k-1,m+k+1}_{1\cdots\ts m,
m+1\cdots\ts m+k}(q^{-l_{m+k,i}})
\een
takes a basis vector $\zeta^{}_{\La}$ of $L(\la)^+_{\mu}$
to another basis vector with a nonzero coef\/f\/icient provided
the respective array $\La+\de_{m+k,i}$ or
$\La-\de_{m+k,i}$ is a pattern.
\epf

\bre\label{rem:restr}
The above proof actually shows that $L(\la)^+_{\mu}$ is irreducible
as a representation of the $q$-Yangian $\Y_q(\gl_n)$.
\par
Moreover, as the coef\/f\/icients of the quantum minors
\ben
{t\ts}^{1\cdots\ts r}_{1\cdots\ts r}(u),\qquad
{t\ts}^{1\cdots\ts r}_{1\cdots\ts r-1,r+1}(u),\qquad
{t\ts}^{1\cdots\ts r-1,r+1}_{1\cdots\ts r}(u)
\een
and
\ben
{\bar t\ts}^{1\cdots\ts r}_{1\cdots\ts r}(u),\qquad
{\bar t\ts}^{1\cdots\ts r}_{1\cdots\ts r-1,r+1}(u),\qquad
{\bar t\ts}^{1\cdots\ts r-1,r+1}_{1\cdots\ts r}(u)
\een
with $r\geqslant 1$ generate the algebra $\U_q(\wh\gl_n)$,
Corollary~\ref{cor:imqmin} together with the formulas
\eqref{qmactgtd}, \eqref{qmactgtr} and \eqref{qmactgtl}
provide an explicit realization of the $\U_q(\wh\gl_n)$-module
$L(\la)^+_{\mu}$.
\ere

Our next goal is to calculate the highest weight and Drinfeld
polynomials of the $\U_q(\wh\gl_n)$-module $L(\la)^+_{\mu}$.
As we have noticed above,
$\zeta_{\La^0}$ is a unique vector
of maximal weight with respect to the standard ordering on the set of weights
of the $\U'_q(\gl_n)$-module $L(\la)^+_{\mu}$.
By the def\/ining relations of $\U_q(\wh\gl_n)$, we have
\ben
t^{\ts(0)}_{jj}\ts t_{ia}(u)
=q^{\de_{ij}-\de_{ja}}\ts t_{ia}(u)\ts t^{\ts(0)}_{jj}
\Fand
t^{\ts(0)}_{jj}\ts \bar t_{ia}(u)=
q^{\de_{ij}-\de_{ja}}\ts\bar t_{ia}(u)\ts t^{\ts(0)}_{jj}.
\een
This implies that $\zeta_{\La^0}$ is the highest vector of the
$\U_q(\wh\gl_n)$-module $L(\la)^+_{\mu}$.
The following theorem provides an identif\/ication of this module;
cf.~\cite{m:yt,nt:ry}. Given three integers $i,j,k$
we shall denote by $\midd\{i,j,k\}$ that of the three which is
between the two others.

\bth\label{thm:hw}
The $\U_q(\wh\gl_n)$-module $L(\la)^+_{\mu}$
is isomorphic to the highest weight representation
$L(\nu(u),\bar\nu(u))$, where the components of the highest weight
are found by
\ben
\nu_k(u)=\frac{(q^{\ts\nu_k^{(1)}}-q^{-\nu_k^{(1)}}\tss u^{-1})\cdots
(q^{\ts\nu_k^{(m+1)}}-q^{-\nu_k^{(m+1)}+2m}\tss u^{-1})}
{(q^{\ts\mu_1}-q^{-\mu_1}\tss u^{-1})\cdots
(q^{\ts\mu_m}-q^{-\mu_m+2m-2}\tss u^{-1})}
\een
and
\ben
\bar\nu_k(u)=\frac{(q^{-\nu_k^{(1)}}-q^{\ts\nu_k^{(1)}}\tss u)\cdots
(q^{-\nu_k^{(m+1)}}-q^{\ts\nu_k^{(m+1)}-2m}\tss u)}
{(q^{-\mu_1}-q^{\ts\mu_1}\tss u)\cdots
(q^{-\mu_m}-q^{\ts\mu_m-2m+2}\tss u)},
\een
where $k=1,\dots,n$ and
\beq\label{nu}
\nu^{(i)}_k=\midd\{\mu_{i-1},\mu_i,\lambda_{k+i-1}\}
\eeq
assuming $\mu_{m+1}=-\infty$, and $\mu_0=+\infty$.
\eth

\bpf
Due to Theorem~\ref{thm:irred},
we only need to calculate the highest weight of
the $\U_q(\wh\gl_n)$-module $L(\la)^+_{\mu}$.
As $\zeta_{\La^0}$ is the highest vector, we have for any
$1\leqslant k\leqslant n$
\ben
{t\ts}^{1\cdots\ts k}_{1\cdots\ts k}(u)\ts \zeta_{\La^0}
=\nu_1(u)\ts \nu_2(q^{-2}\tss u)\cdots \nu_k(q^{-2k+2}\tss u)\ts \zeta_{\La^0}.
\een
On the other hand, by Corollary~\ref{cor:imqmin},
the action of this quantum minor can be found by
applying the evaluation homomorphism \eqref{eval}
to the following series with coef\/f\/icients in $\U_q(\wh\gl_{m+n})$,
\ben
\big(\ts{t\ts}^{1\dots\ts m}_{1\dots\ts m}(u)\big)^{-1}\cdot
{t\ts}^{1\dots\ts m,\ts m+1\dots\ts m+k}_{1\dots\ts m,\ts m+1\dots\ts m+k}
(u),
\een
and then applying the image to the vector $\zeta_{\La^0}$.
This yields the relation
\begin{gather*}
\nu_1(u)\ts \nu_2(q^{-2}\tss u)\cdots \nu_k(q^{-2k+2}\tss u)\\
\qquad{}=\frac{(q^{\ts\la_{m+k,1}^{(0)}}-q^{-\la_{m+k,1}^{(0)}}\tss u^{-1})\cdots
(q^{\ts\la_{m+k,m+k}^{(0)}}-q^{-\la_{m+k,m+k}^{(0)}+2m+2k-2}\tss u^{-1})}
{(q^{\ts\mu_1}-q^{-\mu_1}\tss u^{-1})\cdots
(q^{\ts\mu_m}-q^{-\mu_m+2m-2}\tss u^{-1})}.
\end{gather*}
Hence, replacing here $k$ by $k-1$ and $u$ by $q^{2k-2}\tss u$, we get
\begin{gather*}
\nu_k(u)=\frac{(q^{\ts\la_{m+k,1}^{(0)}}-q^{-\la_{m+k,1}^{(0)}-2k+2}\tss u^{-1})\!\cdots\!
(q^{\ts\la_{m+k,m+k}^{(0)}}-q^{-\la_{m+k,m+k}^{(0)}+2m}\tss u^{-1})}
{(q^{\ts\la_{m+k-1,1}^{(0)}}-q^{-\la_{m+k-1,1}^{(0)}-2k+2}\tss u^{-1})\!\cdots\!
(q^{\ts\la_{m+k-1,m+k-1}^{(0)}}-q^{-\la_{m+k-1,m+k-1}^{(0)}+
2m-2}\tss u^{-1})}.
\end{gather*}
By the def\/inition \eqref{lao} of the pattern $\La^0$, we have
\ben
\la_{m+k,i}^{(0)}=\la_{m+k-1,i}^{(0)}=\la_i\qquad\text{for}\quad i=1,\dots,k-1
\een
while for any $1\leqslant j\leqslant m$ we have
\begin{gather*}
\frac{(q^{\ts\la_{m+k,k+j-1}^{(0)}}-q^{-\la_{m+k,k+j-1}^{(0)}+2j-2}\tss u^{-1})}
{(q^{\ts\la_{m+k-1,k+j-1}^{(0)}}-q^{-\la_{m+k-1,k+j-1}^{(0)}+2j-2}\tss u^{-1})}
\cdot
(q^{\ts\mu_j}-q^{-\mu_j+2j-2}\tss u^{-1})\\
\qquad{}=(q^{\ts\nu_{k}^{(j)}}-q^{-\nu_{k}^{(j)}+2j-2}\tss u^{-1}),
\end{gather*}
completing the proof of the formula for $\nu_k(u)$.

Similarly, using Corollary~\ref{cor:imqmin} again, we f\/ind that
\begin{gather*}
\bar\nu_1(u)\ts \bar\nu_2(q^{-2}\tss u)\cdots \bar\nu_k(q^{-2k+2}\tss u)\\
\qquad{}=\frac{(q^{-\la_{m+k,1}^{(0)}}-q^{\ts\la_{m+k,1}^{(0)}}\tss u)\cdots
(q^{-\la_{m+k,m+k}^{(0)}}-q^{\ts\la_{m+k,m+k}^{(0)}-2m-2k+2}\tss u)}
{(q^{-\mu_1}-q^{\ts\mu_1}\tss u)\cdots
(q^{-\mu_m}-q^{\ts\mu_m-2m+2}\tss u)},
\end{gather*}
and then proceed in the same way as in the calculation of $\nu_k(u)$.
\epf

In the case where the components of $\la$ and $\mu$
are nonnegative we may regard them as partitions.
Consider the corresponding skew diagram $\la/\mu$
and denote by $c(\al)$ the {\it content\/} of a box $\alpha\in\lambda/\mu$
so that $c(\al)=j-i$ if $\al$ occurs in row $i$ and column $j$.

\bco\label{cor:dp}
The Drinfeld polynomials $P_1(u),\dots,P_{n-1}(u)$ corresponding
to the $\U_q(\wh\gl_n)$-module $L(\la)^+_{\mu}$ are given by
\ben
P_k(u)=\prod_{i=1}^{m+1}
(1-q^{2\nu_{k+1}^{(i)}-2i+2}u)(1-q^{2\nu_{k+1}^{(i)}-2i+4}u)\cdots
(1-q^{2\nu_{k}^{(i)}-2i}u),
\een
for $k=1,\dots,n-1$. If $\la$ and $\mu$ are partitions then
the formula can also be written as
\ben
P_k(u)=\prod_{\al} (1-q^{2\ts c(\al)}u),
\een
where $\al$ runs over the top boxes of the columns of height $k$ in
the diagram of $\la/\mu$.
\eco

\bpf
The formulas follow from \eqref{nunupone} and Theorem~\ref{thm:hw}.
\epf

Note that if $m=0$ then $L(\la)^+_{\mu}$ may be regarded
as the evaluation module $L(\la)$ over $\U_q(\wh\gl_n)$.
Clearly, the Drinfeld polynomials provided by Corollary~\ref{cor:dp}
coincide with those given by \eqref{evalmod}.

\bre\label{rem:yang}
The arguments of this section can be applied to produce simpler proofs
of the Yangian analogues of Theorem~\ref{thm:hw}
and Corollary~\ref{cor:dp}; cf.~\cite{m:yt,nt:ry}.
\ere

Finally, we calculate the {\it Gelfand--Tsetlin characters\/}
of the skew representations.
Following \cite[Section~5.2]{bk:rs},
introduce the set $\Pc_n$ of all power series of the form
$a(u)=a_1(u_1)\cdots a_n(u_n)$, where the $u_1,\dots,u_n$ are indeterminates
and each $a_i(u)$ is a power series in $u^{-1}$.
We shall denote by $a_i^{(r)}$ the coef\/f\/icient of $u^{-r}$ of this series.
Consider the group algebra $\ZZ[\Pc_n]$ of the Abelian group~$\Pc_n$
whose elements are f\/inite linear combinations of the form
$\sum m_{a(u)}[a(u)]$, where $m_{a(u)}\in\ZZ$.

We shall be working with the $q$-Yangian $\Y_q(\gl_n)$.
Let us introduce the series $h_i(u)$ with coef\/f\/icients in $\Y_q(\gl_n)$
by the formulas
\ben
h_i(u)={t\ts}^{1\dots \ts i-1}_{1\dots\ts i-1}(q^{2i-2}u)^{-1}\ts
{t\ts}^{1\dots \ts i}_{1\dots\ts i}(q^{2i-2}u),
\qquad i=1,\dots,n.
\een
Due to \eqref{center}, the coef\/f\/icients $h_i^{(r)}$ of all the series
form a commutative subalgebra of $\Y_q(\gl_n)$.
Now, if $V$ is a f\/inite-dimensional representation of $\Y_q(\gl_n)$
and $a(u)\in\Pc_n$, the corresponding Gelfand--Tsetlin subspace
$V_{a(u)}$ consists of the vectors $v\in V$ with the property that for
each $i=1,\dots,n$ and $r\geqslant 0$ there exists $p\geqslant 1$
such that $(h_i^{(r)}-a_i^{(r)})^p\ts v=0$. Then the Gelfand--Tsetlin
character of $V$ is def\/ined by
\ben
\ch V=\sum_{a(u)\in\Pc_n} (\dim V_{a(u)})\tss [a(u)];
\een
cf.~\cite{fr:qc,k:st}. By analogy with \cite[Section~6.2]{bk:rs},
introduce the following special elements of the group algebra $\ZZ[\Pc_n]$
by
\ben
x_{i,a}=\Big[\ts\frac{q^{\tss a+i}-q^{-a-i}\ts u_i^{-1}}
{q^{\tss a+i-1}-q^{-a-i+1}\ts u_i^{-1}}\Big],\qquad 1\leqslant i\leqslant n,
\quad a\in\CC.
\een

We make $L(\la)^+_{\mu}$ into a module over the $q$-Yangian
$\Y_q(\gl_n)$ by taking the composition
of the homomorphism $\Y_q(\gl_n)\to\Y_q(\gl_{m+n})$
given by the f\/irst relation in
\eqref{homnatshar} with the evaluation
homomorphism
$\Y_q(\gl_{m+n})\to \U_q(\gl_{m+n})$
given by
\ben
T(u)\mapsto \frac{T-\overline T\ts u^{-1}}{1-u^{-1}}.
\een
This additional factor, as compared to \eqref{eval}, will ensure the character
formulas look simpler. For the same purpose, we shall also assume that
the components of $\la$ and $\mu$ are non-negative so that
we can consider the corresponding skew diagram
$\la/\mu$. The formulas in the general case will then be
obtained by an obvious modif\/ication.
A {\it semistandard $\la/\mu$-tableau $\Tc$\/} is obtained
by writing the elements of the set $\{1,\dots,n\}$ into
the boxes of the diagram of $\la/\mu$ in such a~way that the
elements in each row weakly increase while the elements in each column
strictly increase. By $\Tc(\al)$ we denote the entry of $\Tc$ in the box
$\al\in\la/\mu$.

\bco\label{cor:character}
The Gelfand--Tsetlin
character of the $\Y_q(\gl_n)$-module $L(\la)^+_{\mu}$
is given by
\ben
\ch L(\la)^+_{\mu}=\sum_{\Tc}\prod_{\al\in\la/\mu} x^{}_{\Tc(\al),c(\al)},
\een
summed over all semistandard $\la/\mu$-tableaux $\Tc$.
\eco

\bpf
First, we consider the particular case $m=0$ so that
$L(\la)^+_{\mu}$ coincides with the evaluation module $L(\la)$.
The coef\/f\/icients of the quantum determinant of $\Y_q(\gl_n)$ act on
$L(\la)$ as scalar operators found from
\ben
\qdet T(u)|^{}_{L(\la)}=\prod_{i=1}^n\ts
\frac{q^{\tss\la_i}-q^{-\la_i+2i-2}u^{-1}}{1-q^{\tss 2i-2}u^{-1}}.
\een
Observe that, regarding $\la$ as a diagram,
we can write the product here as
\ben
\prod_{\al\in\la}
\frac{q^{\tss c(\al)+1}-q^{-c(\al)-1}u^{-1}}{q^{\tss c(\al)}-q^{-c(\al)}u^{-1}}.
\een
Note also that the quantum determinant can be factorized as
\ben
\qdet T(u)=h_1(u)\tss h_2(q^{-2}u)\cdots h_n(q^{-2n+2}u).
\een
Now we employ the well-known bijection between the
patterns associated with $\la$ and the semistandard $\la$-tableaux.
Namely, the pattern $\La$ can be viewed as the sequence of diagrams
\ben
\La_1\subseteq \La_2\subseteq\cdots
\subseteq \La_n=\lambda,
\een
where the $\La_k$ are the rows of $\La$.
The betweenness conditions \eqref{betw} mean that the skew diagram
$\La_k/\La_{k-1}$ is a {\it horizontal strip\/};
see, e.g., Macdonald~\cite[Chapter~1]{m:sf}.
The corresponding semistandard tableau is obtained by placing
the entry $k$ into each box of $\La_k/\La_{k-1}$.
By the above observation,
for any basis vector $\xiL\in L(\la)$ and any $1\leqslant k\leqslant n$
we have
\ben
h_1(u)\tss h_2(q^{-2}u)\cdots h_k(q^{-2k+2}u)\ts\xiL
=\prod_{\al\in\La_k}
\frac{q^{\tss c(\al)+1}-q^{-c(\al)-1}u^{-1}}{q^{\tss c(\al)}-q^{-c(\al)}u^{-1}}
\ts\xiL,
\een
so that
\beql{hpract}
h_k(q^{-2k+2}u)\ts\xiL=\prod_{\al\in\La_k/\La_{k-1}}
\frac{q^{\tss c(\al)+1}-q^{-c(\al)-1}u^{-1}}{q^{\tss c(\al)}-q^{-c(\al)}u^{-1}}
\ts\xiL.
\eeq
Thus, the element of $\ZZ[\Pc_n]$ corresponding to the
eigenvalue of $h_k(u_k)$ on $\xiL$ coincides with the product
\ben
\prod_{\al\in\La_k/\La_{k-1}} x^{}_{k,c(\al)}.
\een
This proves the claim for the case $m=0$.

In the skew case ($m\geqslant 1$)
we use Corollary~\ref{cor:imqmin} which implies that
formula \eqref{hpract} remains valid when the action is considered
on the basis vector $\zeta^{}_{\La}\in L(\la)^+_{\mu}$
instead of $\xiL$. Here we
denote by $\La_k$ the row $(\la_{m+k,1},\dots,\la_{m+k,m+k})$
of the trapezium pattern $\La$ and
use a natural bijection
between the trapezium patterns and the semistandard $\la/\mu$-tableaux.
\epf

\bre\label{rem:yangco}
The Gelfand--Tsetlin character of the
corresponding skew module $L(\la)^+_{\mu}$
over the Yangian $\Y(\gl_n)$ is given by the same formula as in
Corollary~\ref{cor:character}, where $x_{i,a}$ is now def\/ined
by $x_{i,a}=[1+(u_i+a+i-1)^{-1}]$; see \cite[Section~6.2]{bk:rs}.
The skew representations of the Yangians and quantum af\/f\/ine algebras
were studied in the literature
from various viewpoints providing dif\/ferent interpretations
of the character formula of Corollary~\ref{cor:character}; see
e.g.~\cite{br:rs,c:ni,m:yt,nn:pt,nt:ry}. For the evaluation modules,
a calculation similar to the above
can be found e.g. in
\cite[Section~4.5]{fm:ha}, \cite[Section~7.4]{bk:rs}.

One can easily extend the def\/inition of the Gelfand--Tsetlin character
to representations of $\U_q(\wh\gl_n)$ by considering
the commutative subalgebra generated by the coef\/f\/icients of the
series~$h_i(u)$ together with the coef\/f\/icients of the $\bar h_i(u)$ which
are def\/ined in the same way with the use of the quantum minors of
$\overline T(u)$. However, by the results of \cite{fr:qc}, the eigenvalues
of the $\bar h_i(u)$ on f\/inite-dimensional
representations are essentially determined by those of $h_i(u)$.
In particular, the character formula for the $\U_q(\wh\gl_n)$-module
$L(\la)^+_{\mu}$ will have the form given in Corollary~\ref{cor:character}.
\ere

\section[Olshanski algebra associated with $\U_q(\gl_n)$]{Olshanski algebra associated
with $\boldsymbol{\U_q(\gl_n)}$}\label{sec:oa}
\setcounter{equation}{0}

In this section we consider $q$ as a formal variable so that
the quantum algebras are regarded as algebras over $\CC(q)$.
Denote by $\wt\U_q(\gl_n)$ the subalgebra of $\U_q(\gl_n)$
generated by the elements
\beql{tauint}
\tau_{ij}=t_{ij}\tss \bar t_{jj},\quad i>j\Fand
\bar\tau_{ij}=\bar t_{ij}\tss \bar t_{jj},\quad i\leqslant j.
\eeq
Since in the algebra $\U_q(\gl_n)$ we have
\ben
q^{\de_{ja}}\ts t_{ia}\ts \bar t_{jj}=
q^{\de_{ij}}\ts \bar t_{jj}\ts t_{ia}
\Fand
q^{\de_{ja}}\ts \bar t_{ia}\ts \bar t_{jj}=
q^{\de_{ij}}\ts \bar t_{jj}\ts \bar t_{ia},
\een
we may regard $\wt\U_q(\gl_n)$ as an associative algebra
generated by the elements $\tau_{ij}$ with $i>j$
and~$\bar\tau_{ij}$ with $i\leqslant j$ subject to the def\/ining
relations
\begin{gather*}
q^{\de_{ij}+\de_{ja}}\ts \tau_{ia}\ts \tau_{jb}-
q^{\de_{ib}+\de_{ab}}\ts \tau_{jb}\ts \tau_{ia}
=(q-\qin)\ts q^{\de_{ia}} (\de_{b<a} -\de_{i<j})
\ts \tau_{ja}\ts \tau_{ib},\\
q^{\de_{ij}+\de_{ja}}\ts \bar\tau_{ia}\ts \bar\tau_{jb}-
q^{\de_{ib}+\de_{ab}}\ts \bar\tau_{jb}\ts \bar\tau_{ia}
=(q-\qin)\ts q^{\de_{ia}} (\de_{b<a} -\de_{i<j})
\ts \bar\tau_{ja}\ts \bar\tau_{ib},\\
q^{\de_{ij}+\de_{ja}}\ts \bar \tau_{ia}\ts \tau_{jb}-
q^{\de_{ib}+\de_{ab}}\ts\tau_{jb}\ts \bar \tau_{ia}
=(q-\qin)\ts q^{\de_{ia}} \ts \big(\de_{b<a}\ts
\tau_{ja}\ts \bar\tau_{ib} -\de_{i<j}\ts
\ts \bar\tau_{ja}\ts \tau_{ib}\big),
\end{gather*}
where
\ben
\tau_{ij}=\bar \tau_{ji}=0, \qquad 1 \leqslant i<j\leqslant n,
\Fand \tau_{ii}=1,\qquad i=1,\dots,n.
\een

For any positive integer $n$
the algebra $\wt\U_q(\gl_{n-1})$ can be identif\/ied with a subalgebra
of $\wt\U_q(\gl_n)$ generated by the elements $\tau_{ij}$ and $\bar\tau_{ij}$ with
$1\leqslant i,j\leqslant n-1$.

Fix a nonnegative integer $m$ such that $m\leqslant n$ and
denote by $\wt\U_q(\gl_{n,m})$ the subalgebra of $\wt\U_q(\gl_n)$
generated by the elements $\tau_{ij}$ and $\bar\tau_{ij}$ with
$m+1\leqslant i,j\leqslant n$. This subalgebra is isomorphic to
$\wt\U_q(\gl_{n-m})$. Let $\Ar_m(n)$ denote the centralizer
of $\wt\U_q(\gl_{n,m})$ in $\wt\U_q(\gl_n)$.
Also, let $\Ar(n)^0$ denote the centralizer of the element $\bar\tau_{nn}$
in $\wt\U_q(\gl_n)$ and let $\Ir(n)$ be the left ideal in $\wt\U_q(\gl_n)$
generated by the
elements $\bar\tau_{in}$, $i=1,\dots,n$. Then the Poincar\'e--Birkhof\/f--Witt
theorem for the algebra $\U_q(\gl_n)$ implies that
$\Ir(n)^0=\Ir(n)\cap \Ar(n)^0$
is a two-sided ideal in $\Ar(n)^0$ and one has a~vector space
decomposition
\ben
\Ar(n)^0=\Ir(n)^0\oplus \wt\U_q(\gl_{n-1}).
\een
Therefore, the projection of $\Ar(n)^0$ onto $\wt\U_q(\gl_{n-1})$
with the kernel $\Ir(n)^0$ is an algebra homomorphism.
If $m<n$ then its
restriction to the subalgebra $\Ar_m(n)$ def\/ines
a homomorphism
\begin{gather}\label{proj}
o_n: \Ar_m(n)\to \Ar_m(n-1).
\end{gather}
Note that the algebra $\Ar_m(n)$ inherits the f\/iltration of $\wt\U_q(\gl_n)$
def\/ined by
\ben
\deg \tau_{ij}=0,\quad i>j,\Fand \deg \bar\tau_{ij}=1,\quad i\leqslant j.
\een
Clearly, the homomorphism $o_n$ is f\/iltration preserving.

\bde\label{def:oam}
The {\it Olshanski algebra\/} $\Ar_m$ is defined as the
projective limit of the sequence of the algebras $\Ar_m(n)$, $n\geqslant m$,
with respect to the homomorphisms
\ben
\begin{CD}
\Ar_m(m)
@<{o_{m+1}}<<
\Ar_m(m+1)
@<{o_{m+2}}<<\cdots
@<{o_{n}}<<
\Ar_m(n)@<{o_{n+1}}<<\cdots,
\end{CD}
\een
where the limit
is taken in the category of filtered associative algebras.
\ede

An element of the algebra $\Ar_m$ is a sequence
of the form
$a=(a_m,a_{m+1},\dots,a_n,\dots)$ with $a_n\in\Ar_m(n)$,
$o_{n}(a_n)=a_{n-1}$ for $n>m$, and
\beql{degina}
\deg a=\sup_{n\geqslant m}\deg a_n<\infty,
\eeq
where
$\deg a_n$ denotes the degree of $a_n$ in $\wt\U_q(\gl_n)$.
If $b=(b_m,b_{m+1},\dots,b_n,\dots)$ is another element of $\Ar_m$
then the product $ab$ is the sequence
\ben
ab=(a_mb_m,a_{m+1}b_{m+1},\dots,a_nb_n,\dots).
\een

We def\/ine the algebra $\wt\U_q(\gl_{\infty})$ as the inductive limit
of the algebras $\wt\U_q(\gl_n)$ with respect to the natural embeddings
$\wt\U_q(\gl_n)\hra\wt\U_q(\gl_{n+1})$,
\ben
\wt\U_q(\gl_{\infty})=\underset{n\geqslant 1}{\bigcup} \ts\wt\U_q(\gl_n).
\een
Note that the algebra $\Ar_0(n)$ coincides with the center
of $\wt\U_q(\gl_n)$.
The elements of the algebra~$\Ar_0$ can therefore be regarded as
{\it virtual Casimir elements\/} for the algebra $\wt\U_q(\gl_{\infty})$.

Let us apply the evaluation homomorphism
\beql{evalqy}
\Y_q(\gl_n)\to\U_q(\gl_n),\qquad
T(u)\mapsto T-\overline T\tss u^{-1},
\eeq
see \eqref{eval}, to the quantum determinant $\qdet T(u)$. We get
the polynomial in $u^{-1}$,
\begin{gather}
\qdet(T-\overline T\tss u^{-1})
=\sum_{\si\in\Sym_n}(-q)^{-l(\si)}\tss
(t_{\si(1)1}-\bar t_{\si(1)1}\tss u^{-1})\nonumber\\
\phantom{\qdet(T-\overline T\tss u^{-1})=}{}\times(t_{\si(2)2}-\bar t_{\si(2)2}\tss q^2\tss u^{-1})\cdots
(t_{\si(n)n}-\bar t_{\si(n)n}\tss q^{2n-2}\tss u^{-1}),\label{qdettt}
\end{gather}
whose coef\/f\/icients are central in $\U_q(\gl_n)$.
Applying the Harish-Chandra homomorphism \eqref{hchhom}
to the coef\/f\/icients (see Section~\ref{sec:def}), we f\/ind
\ben
\chi:\qdet(T-\overline T\tss u^{-1})\mapsto q^{n(n-1)/2}\cdot
(x_1-x_1^{-1}u^{-1})\cdots (x_n-x_n^{-1}u^{-1}).
\een
Hence, the coef\/f\/icients of the polynomial \eqref{qdettt}
generate the center $\Zr_q$ of $\U_q(\gl_n)$.

The product
\ben
d_n(u)=\qdet(T-\overline T\tss u^{-1})\cdot\bar t_{11}\cdots \bar t_{nn}
\een
is a polynomial in $u^{-1}$ with constant term $1$ whose coef\/f\/icients belong to
the center of the subalgebra $\wt\U_q(\gl_n)$. Explicitly, this polynomial
can be written as
\beql{dn}
d_n(u)
=\sum_{\si\in\Sym_n}(-q)^{-l(\si)}\tss q^{\text{ind}(\si)}\tss
(\tau_{\si(1)1}-\bar \tau_{\si(1)1}\tss u^{-1})
\cdots
(\tau_{\si(n)n}-\bar \tau_{\si(n)n}\tss q^{2n-2}\tss u^{-1}),
\eeq
where
$
\text{ind}(\si)=\sharp\ts\{i=1,\dots,n\ |\ \si(i)<i\}.
$
Write
\ben
d_n(u)=1+d_n^{\tss(1)}\tss u^{-1}+\dots+d_n^{\tss(n)}\tss u^{-n}.
\een

\bpr\label{prop:gencen}
The elements $d_n^{\tss(1)},\dots,d_n^{\tss(n)}$
are algebraically independent and generate
the center of the algebra $\wt\U_q(\gl_n)$.
\epr

\bpf
Any central element $z$ of $\wt\U_q(\gl_n)$ must commute with
$\bar\tau_{11},\dots,\bar\tau_{nn}$. Using the isomorphism \eqref{tenprudec},
and writing $z$ as a linear combination of the corresponding basis monomials,
we conclude that $z$ must also commute with $t_{11},\dots,t_{nn}$.
Therefore, $z$ belongs to the center $\Zr_q$ of~$\U_q(\gl_n)$.
Hence, $z$ is a polynomial in $t^{}_0,t_0^{-1},d_n^{\tss(1)},\dots,d_n^{\tss(n)}$,
where $t^{}_0=t_{11}\cdots t_{nn}$. However, such a polynomial
does not belong to the subalgebra $\wt\U_q(\gl_n)$ unless
it only contains non-positive even powers of $t^{}_0$.
Since $t^{-2}_0$ coincides with $d_n^{\tss(n)}$ up to a nonzero
constant factor, this proves that the center of $\wt\U_q(\gl_n)$
is generated by $d_n^{\tss(1)},\dots,d_n^{\tss(n)}$.

Finally, the image of $d_n(u)$ under the Harish-Chandra homomorphism
is given by
\ben
\chi: d_n(u)\mapsto (1-x_1^{-2}\tss u^{-1})\cdots (1-x_n^{-2}\tss u^{-1}),
\een
thus proving the algebraic independence of $d_n^{\tss(1)},\dots,d_n^{\tss(n)}$.
\epf

Applying the homomorphism \eqref{proj} with $m=0$
to the coef\/f\/icients of
the polynomial $d_n(u)$ we f\/ind immediately that
\ben
o_n:d_n(u)\mapsto d_{n-1}(u).
\een
Since the degree of the element $d_n^{\tss(k)}$ does not exceed $k$
for all $n\geqslant k$, we may def\/ine a virtual
Casimir element $d^{\tss(k)}\in\Ar_0$ for any $k\geqslant 1$
as the sequence
\ben
d^{\tss(k)}=(d_n^{\tss(k)}\ |\ n\geqslant 0),
\een
where we set $d_n^{\tss(k)}=0$ for $n<k$.
In order to get an alternative expression for
the virtual Casimir elements $d^{\tss(k)}$,
denote by $\Sym_{\infty}$ the group of f\/inite permutations
of the set of positive integers, so that for any
$p\in\Sym_{\infty}$ we have $p(l)=l$ for all suf\/f\/iciently large $l$.
Def\/ine the {\it virtual quantum determinant\/} as the
formal power series
\ben
d(u)=1+d^{\tss(1)}u^{-1}+d^{\tss(2)}u^{-2}+\cdots.
\een
Using the expression
\ben
d(u)
=\sum_{\si\in\Sym_{\infty}}(-q)^{-l(\si)}\tss q^{\text{ind}(\si)}\tss
(\tau_{\si(1)1}-\bar \tau_{\si(1)1}\tss u^{-1})
(\tau_{\si(2)2}-\bar \tau_{\si(2)2}\tss q^{2}\tss u^{-1})\cdots,
\een
we can regard the coef\/f\/icients $d^{\tss(k)}$ as certain formal series
of elements of $\wt\U_q(\gl_{\infty})$.

The following description of the algebra $\Ar_0$
is implied by Proposition~\ref{prop:gencen}.

\bpr\label{prop:vqd}
The elements $d^{\tss(1)},d^{\tss(2)},\dots$ are algebraically independent
and generate the algebra $\Ar_0$.
\epr

We are now in a position to establish a relationship between
the Olshanski algebra $\Ar_m$
and the $q$-Yangian $\Y_q(\gl_m)$.
It will be convenient to work with the subalgebra $\wt\Y_q(\gl_m)$
of $\Y_q(\gl_m)$ generated by the coef\/f\/icients of the series
\beql{taudef}
\tau_{ij}(u)=t_{ij}(u)\ts \bar t_{jj}^{\ts(0)},\qquad 1\leqslant i,j\leqslant m.
\eeq
As with the algebra $\wt\U_q(\gl_m)$,
using the relations
\beql{reztze}
q^{\de_{ja}}\ts t_{ia}(u)\ts \bar t_{jj}^{\ts(0)}=
q^{\de_{ij}}\ts \bar t_{jj}^{\ts(0)}\ts t_{ia}(u)
\Fand
q^{\de_{ja}}\ts \bar t_{ia}(u)\ts \bar t_{jj}^{\ts(0)}=
q^{\de_{ij}}\ts \bar t_{jj}^{\ts(0)}\ts \bar t_{ia}(u),
\eeq
it is easy to write down the def\/ining relations of $\wt\Y_q(\gl_m)$
in terms of the coef\/f\/icients $\tau_{ij}^{(r)}$ of the series $\tau_{ij}(u)$.

Now we use the quantum Sylvester theorem (Theorem~\ref{thm:qs}).
Taking the composition of the homomorphism \eqref{homnatut}
and the evaluation homomorphism \eqref{evalqy},
we obtain a homomorphism $\Y_q(\gl_m)\to\U_q(\gl_n)$ which can be
written as
\beql{tijmin}
t_{ij}(u)\mapsto
\qdet(T-\overline T\tss u^{-1})^{}_{\Bc_i\Bc_j},
\eeq
where we use the quantum determinant
of the submatrix of $T-\overline T\tss u^{-1}$ corresponding to the
rows~$\Bc_i$ and columns $\Bc_j$,
where $\Bc_i$ denotes the set $\{i,m+1,\dots,n\}$.
By the relations \eqref{center},
the image of the series $t_{ij}(u)$ under the homomorphism
\eqref{tijmin} commutes with the elements $t_{kl}$ and~$\bar t_{kl}$
with $m+1\leqslant k,l\leqslant n$.
Introduce the matrices $\Tc=[\tau_{ij}]$ and $\overline\Tc=[\bar\tau_{ij}]$
with $1\leqslant i,j\leqslant n$ and set
\ben
\qdet'(\Tc-\overline \Tc\tss u^{-1})^{}_{\Bc_i\Bc_j}=
\qdet(T-\overline T\tss u^{-1})^{}_{\Bc_i\Bc_j}
\cdot \bar t_{jj}\ts\bar t_{m+1,m+1}\cdots \bar t_{nn}.
\een
An explicit formula for this quantum minor
in terms of $\tau_{ij}$ and $\bar\tau_{ij}$ can be written
with the use of the parameter $\text{ind}(\si)$ as in \eqref{dn}.
Since the product $\bar t_{m+1,m+1}\dots \bar t_{nn}$
commutes with $\qdet(T-\overline T\tss u^{-1})^{}_{\Bc_i\Bc_j}
\cdot \bar t_{jj}$ (see \eqref{center}),
the mapping
\beql{wtyu}
\vp_n:\tau_{ij}(u)\mapsto
\qdet'(\Tc-\overline \Tc\tss u^{-1})^{}_{\Bc_i\Bc_j}
\eeq
def\/ines a homomorphism $\wt\Y_q(\gl_m)\to\wt\U_q(\gl_n)$.
Furthermore, the same relation \eqref{center} implies that
the image of the homomorphism \eqref{wtyu} is contained in the centralizer
$\Ar_m(n)$ so that we have a homomorphism
\ben
\vp_n:\wt\Y_q(\gl_m)\to\Ar_m(n).
\een
Using the explicit formula for
$\qdet'(\Tc-\overline \Tc\tss u^{-1})^{}_{\Bc_i\Bc_j}$
as in \eqref{dn}, we f\/ind that the diagram
\ben
\begin{CD}
\wt\Y_q(\gl_m) @= \wt\Y_q(\gl_m) @= \cdots @= \wt\Y_q(\gl_m)@=\cdots\\
@V \vp_m VV     @V\vp_{m+1}VV @.
@V\vp_{n}VV\\
\Ar_m(m) @<<o_{m+1}< \Ar_{m}(m+1)@<<<\cdots
@<<o_n< \Ar_m(n)@<<o_{n+1}<\cdots
\end{CD}
\een
is commutative. Note that the image of the generator $\tau_{ij}^{(r)}$
of $\wt\Y_q(\gl_m)$
under $\vp_n$ has degree ${}\leqslant r$ for any $n$.
Hence the sequence of homomorphisms $(\vp_n\ |\ n\geqslant m)$
def\/ines an algebra homomorphism
$\vp:\wt\Y_q(\gl_m)\to\Ar_m$ which can be written in terms
of the virtual quantum determinants by
\ben
\vp: \tau_{ij}(u)\mapsto \qdet'(\Tc-\overline \Tc\tss u^{-1})^{}_{\Bc_i\Bc_j},
\een
where $\Bc_i$ now denotes the inf\/inite set $\{i,m+1,m+2,\dots\}$.

\bth\label{thm:secoadet}
The homomorphism $\vp$ is an algebra embedding
of $\wt\Y_q(\gl_m)$ into the algebra $\Ar_m$.
\eth

\bpf
We shall use a weak form of
the Poincar\'e--Birkhof\/f--Witt theorem for the algebra
$\wt\Y_q(\gl_m)$ which can be verif\/ied by a direct argument;
see e.g. \cite[Lemma~3.2]{mrs:cs} where a similar result
is proved for a twisted version of the $q$-Yangian.
By this theorem, the ordered monomials in the generators $\tau_{ij}^{(r)}$
span the algebra $\wt\Y_q(\gl_m)$.
It will be suf\/f\/icient to prove that the images of these monomials
under the homomorphism $\vp$ are linearly independent.

As in Section~\ref{sec:def}, set $\A=\CC[q,\qin]$ and consider
the $\A$-subalgebra $\wt\U_{\A}$ of $\wt\U_q(\gl_n)$ generated by
the elements $\tau_{ij}$ and $\bar\tau_{ij}$.
Denote by $\Pc_n$ the algebra of polynomials
in independent variables $x_{ij}$ with
$1\leqslant i,j\leqslant n$.
By the Poincar\'e--Birkhof\/f--Witt theorem for the
algebra $\U_q(\gl_n)$, we have an isomorphism
\beql{limimy}
\wt\U_{\A}\ot_{\A}\CC\cong\Pc_n
\eeq
given by
\ben
\tau_{ij}\mapsto x_{ij}\quad\text{for}\quad
i>j\Fand\bar\tau_{ij}\mapsto x_{ij}\quad\text{for}\quad
i\leqslant j,
\een
where the action of
$\A$ on $\CC$ is def\/ined via the evaluation $q=1$.

Suppose now that the image under $\vp$ of a nontrivial linear combination of
the ordered monomials in the generators $\tau_{ij}^{(r)}$
is zero. We may assume that all coef\/f\/icients of this
linear combination belong to $\A$. Moreover, we may also assume
that at least one coef\/f\/icient does not vanish at $q=1$
so that the image of this linear combination under
the isomorphism \eqref{limimy} yields a nontrivial linear combination of
the corresponding elements of $\Pc_n$.
Observe that the image of the
quantum determinant $\qdet'(\Tc-\overline \Tc\tss u^{-1})^{}_{\Bc_i\Bc_j}$
occurring in \eqref{wtyu} under the isomorphism \eqref{limimy}
coincides with the usual determinant
$\det\tss(X-\overline X\tss u^{-1})^{}_{\Bc_i\Bc_j}$, where
$X$ denote the lower triangular matrix with entries $x_{ij}$ for $i>j$
and with all diagonal entries equal to $1$ while $\overline X$
is the upper triangular matrix with entries $x_{ij}$ for $i\leqslant j$.
Write
\ben
\det\tss(X-\overline X\tss u^{-1})^{}_{\Bc_i\Bc_j}
=\la^{(0)}_{ij}-\la^{(1)}_{ij}\tss u^{-1}+\dots+
(-1)^{n-m+1}\la^{(n-m+1)}_{ij}\tss u^{-n+m-1}.
\een
The proof will be completed
if we show that given any positive integer $p$,
there exists a suf\/f\/iciently large
value of $n$ such that the polynomials
$\la^{(r)}_{ij}$ with $0\leqslant r\leqslant p$ and $1\leqslant i,j\leqslant m$
(assuming $i>j$ for $r=0$) are algebraically independent.
This will lead to a contradiction with the assumption that the ordered
monomials in the $\tau_{ij}^{(r)}$ are linearly dependent.

We shall be proving that the corresponding map
\beql{mapla}
\Lambda_n:\CC^{n^2}\to\CC^{p\tss m^2+m(m-1)/2},
\eeq
given by
\ben
(\tss x_{ij}\ |\ 1\leqslant i,j\leqslant n)\mapsto
(\tss \la^{(r)}_{kl}\ |\ 1\leqslant k,l\leqslant m,\quad 0\leqslant r\leqslant p)
\een
with $k>l$ for $r=0$, is surjective.
We need two auxiliary lemmas.
For any $r\geqslant 1$ we let
$e_r(z_1,\dots,z_l)$ denote the $r$-th elementary
symmetric function in variables $z_1,\dots,z_l$ with $l\geqslant r$
so that
\ben
e_r(z_1,\dots,z_l)=\sum_{1\leqslant i_1<\dots<i_r\leqslant l}
z_{i_1}\cdots z_{i_r}.
\een
We also set $e_0(z_1,\dots,z_l)=1$ and $e_r(z_1,\dots,z_l)=0$ for $r<0$.
Let $l$ be a positive integer and
let $\al_1,\dots,\al_l$ be distinct complex parameters.
For each $i=1,\dots,l$ set
\ben
e_{ri}=e_r(\al_1,\dots,\wh\al_i,\dots,\al_l),
\een
where the hat
indicates the symbol to be omitted. For any $1\leqslant p\leqslant l$
consider the matrix
\ben
\Ec=\begin{pmatrix}
e_{01}&\cdots&e_{0l}\\
e_{11}&\cdots&e_{1l}\\
\vdots&\ddots&\vdots\\
e_{p-1,1}&\cdots&e_{p-1,l}
\end{pmatrix}
\een

\ble\label{lem:minorsnz}
All $p\times p$ minors of the matrix $\Ec$ are nonzero.
\ele

\bpf
Introduce the polynomials $Q_1(t),\dots,Q_l(t)$ by
\ben
Q_i(t)=e_{0i}\tss t^{\tss l-1}+e_{1i}\tss t^{\tss l-2}+\dots+e_{l-1,i}=
\prod_{k=1,\ts k\ne i}^l(t+\al_k).
\een
They are linearly independent for
if $c_1\tss Q_1(t)+\dots+c_l\tss Q_l(t)=0$ then taking $t=-\al_i$
gives $c_i=0$ since the parameters $\al_i$ are all distinct.

Suppose now that columns $i_1,\dots,i_p$ with
$i_1<\dots<i_p$ of the matrix $\Ec$
are linearly dependent. This implies that a certain nontrivial
linear combination
$Q(t)=d_1\tss Q_{i_1}(t)+\dots+d_p\tss Q_{i_p}(t)$ is a~polynomial in $t$
of degree not exceeding $l-p-1$. However, $Q(t)$ has $l-p$ distinct roots
$t=-\al_j$ with $j\in\{1,\dots,l\}\setminus\{i_1,\dots,i_p\}$, and so $Q(t)=0$.
This contradicts the linear independence
of the polynomials $Q_i(t)$.
\epf

Suppose now that $\Ec_0,\dots,\Ec_{m-1}$ are given non-singular
$p\times p$ matrices with complex entries.

\ble\label{lem:vanderm}
For any distinct complex numbers $\be_1,\dots,\be_m$
the block matrix
\beql{blockma}
\begin{pmatrix}\ts
\Ec_0&\be_1\tss\Ec_1&\cdots&\be_1^{m-1}\tss\Ec_{m-1}\\
\Ec_0&\be_2\tss\Ec_1&\cdots&\be_2^{m-1}\tss\Ec_{m-1}\\
\vdots&\vdots&\ddots&\vdots\\
\Ec_0&\be_m\tss\Ec_1&\cdots&\be_m^{m-1}\tss\Ec_{m-1}
\end{pmatrix}
\eeq
is non-singular.
\ele

\bpf
Suppose that a certain linear combination of
the rows of the matrix with coef\/f\/icients $c_1,c_2,\dots,c_{pm}$ is zero.
Since the rows of each of the matrices $\Ec_j$ are linearly
independent, for any index $1\leqslant i\leqslant p$ we
have the linear relations
\ben
\be_1^{\tss k}\tss c_i+\be_2^{\tss k}\tss c_{p+i}+\dots+
\be_m^{\tss k}\tss c_{(m-1)\tss p+i}=0,
\qquad k=0,1,\dots,m-1.
\een
As the $\be_j$ are distinct, for each $i$ the system has only
the trivial solution $c_i=c_{p+i}=\dots=c_{(m-1)\tss p+i}=0$.
\epf

Now we shall show that the map \eqref{mapla} is surjective
for any $n\geqslant (p+1)m$. In fact, we show that this map
is surjective even when the variables $x_{kl}$ with
$m+1\leqslant k\leqslant n$
are specialized
in the following way:
\ben
x_{kl}=\beta_{kl},\qquad\text{for}\quad l=1,\dots,m
\een
and
\ben
x_{kl}=\de_{kl}\al_{k},\qquad\text{for}\quad l=m+1,\dots,n
\een
where $\al_{m+1},\dots,\al_n$ are distinct complex numbers
and the $\be_{kl}$ are certain complex numbers to be chosen below.
Under this specialization, for each f\/ixed value of
the index $i\in\{1,\dots,m\}$
consider the following restriction of the map \eqref{mapla}
\beql{maplai}
\Lambda^{(i)}_n:\CC^n\to\CC^{p\tss m+i-1},
\eeq
given by
\begin{gather*}
(x_{i1},x_{i2},\dots,x_{in})\mapsto
(\tss \la^{(0)}_{i1},\dots,\la^{(p)}_{i1},\dots,
\la^{(0)}_{i,i-1},\dots,\la^{(p)}_{i,i-1},\la^{(1)}_{ii},\dots,\la^{(p)}_{ii},
\dots,\la^{(1)}_{i\tss m},\dots,\la^{(p)}_{i\tss m}),
\end{gather*}
so that the functions $\la^{(r)}_{ij}$ are def\/ined by the
expansions
\ben
\left|\begin{matrix}tx_{ij}&x_{i,m+1}&x_{i,m+2}&\cdots&x_{in}\\
t\be_{m+1,j}&t+\al_{m+1}&0&\cdots&0\\
t\be_{m+2,j}&0&t+\al_{m+2}&\cdots&0\\
\vdots&\vdots&\vdots&\ddots&\vdots\\
t\be_{n,j}&0&0&\cdots&t+\al_n
\end{matrix}\right|
=t^{n-m+1}\la^{(0)}_{ij}+\dots+
t\la^{(n-m)}_{ij},
\een
for $j=1,\dots,i-1$,
\ben
\left|\begin{matrix}t+x_{ii}&x_{i,m+1}&x_{i,m+2}&\cdots&x_{in}\\
t\be_{m+1,i}&t+\al_{m+1}&0&\cdots&0\\
t\be_{m+2,i}&0&t+\al_{m+2}&\cdots&0\\
\vdots&\vdots&\vdots&\ddots&\vdots\\
t\be_{n,i}&0&0&\cdots&t+\al_n
\end{matrix}\right|
=t^{n-m+1}+t^{n-m}\la^{(1)}_{ii}+\dots+
\la^{(n-m+1)}_{ii},
\een
and
\ben
\left|\begin{matrix}x_{ij}&x_{i,m+1}&x_{i,m+2}&\cdots&x_{in}\\
t\be_{m+1,j}&t+\al_{m+1}&0&\cdots&0\\
t\be_{m+2,j}&0&t+\al_{m+2}&\cdots&0\\
\vdots&\vdots&\vdots&\ddots&\vdots\\
t\be_{n,j}&0&0&\cdots&t+\al_n
\end{matrix}\right|
=t^{n-m}\la^{(1)}_{ij}+\dots+
\la^{(n-m+1)}_{ij},
\een
for $j=i+1,\dots,m$. More explicitly, we have
\ben
\la^{(r)}_{ij}=x_{ij}\tss e_{r}(\al_{m+1},\dots,\al_n)-
\sum_{k=m+1}^n\be_{kj}\tss x_{ik}\tss e_{r-1}(\al_{m+1},\dots,\wh\al_k,\dots,\al_n)
\een
for $j=1,\dots,i-1$ and $r\geqslant 0$,
\begin{gather*}
\la^{(r)}_{ii}=x_{ii}\tss e_{r-1}(\al_{m+1},\dots,\al_n)\\
\phantom{\la^{(r)}_{ii}=}{}-\sum_{k=m+1}^n\be_{ki}\tss x_{ik}\tss e_{r-1}(\al_{m+1},\dots,\wh\al_k,\dots,\al_n)
+e_{r}(\al_{m+1},\dots,\al_n)
\end{gather*}
for $r\geqslant 1$, and
\ben
\la^{(r)}_{ij}=x_{ij}\tss e_{r-1}(\al_{m+1},\dots,\al_n)-
\sum_{k=m+1}^n\be_{kj}\tss x_{ik}\tss e_{r-1}(\al_{m+1},\dots,\wh\al_k,\dots,\al_n)
\een
for $j=i+1,\dots,m$ and $r\geqslant 1$.
Since the $\la^{(r)}_{ij}$ are linear functions
in the variables $x_{i1},\dots,x_{in}$,
in order to establish
the surjectivity of the map \eqref{maplai},
it will be suf\/f\/icient
to demonstrate that the rank of the corresponding coef\/f\/icient
matrix is maximal,
that is, equal to $p\tss m+i-1$. Writing down the matrix in an explicit
form and using the observation that $\la^{(0)}_{ij}=x_{ij}$
for $j=1,\dots,i-1$ we conclude that the claim will follow
by proving that the matrix
\beql{imatr}
\begin{pmatrix}
-\be_{m+1,1}e_{0,m+1}&\cdots&-\be_{n1}e_{0,n}\\
-\be_{m+1,1}e_{1,m+1}&\cdots&-\be_{n1}e_{1,n}\\
\vdots&\cdots&\vdots\\
-\be_{m+1,1}e_{p-1,m+1}&\cdots&-\be_{n1}e_{p-1,n}\\
\vdots&\cdots&\vdots\\
\vdots&\cdots&\vdots\\
-\be_{m+1,m}e_{0,m+1}&\cdots&-\be_{nm}e_{0,n}\\
-\be_{m+1,m}e_{1,m+1}&\cdots&-\be_{nm}e_{1,n}\\
\vdots&\cdots&\vdots\\
-\be_{m+1,m}e_{p-1,m+1}&\cdots&-\be_{nm}e_{p-1,n}
\end{pmatrix}
\eeq
has rank $p\tss m$, where $e_{rk}=e_r(\al_{m+1},\dots,\wh\al_k,\dots,\al_n)$.
Now, we use the assumption $n-m\geqslant p\tss m$ and Lemma~\ref{lem:minorsnz}.
Giving appropriate values to the parameters $\be_{kl}$, we may
choose a submatrix of \eqref{imatr} of size $p\tss m\times p\tss m$
of the form \eqref{blockma}, where all the $p\times p$ blocks
$\Ec_i$ are non-singular and the parameters $\be_i$ are distinct.
Thus, the claim follows by the application of Lemma~\ref{lem:vanderm}.

Finally, by allowing the index $i$ in
\eqref{maplai} to vary and by choosing the parameters
$\beta_{kl}$ in the same way for each value of $i\in\{1,\dots,m\}$,
we get a surjective map
\ben
\CC^{m\tss n}\to\CC^{p\tss m^2+m(m-1)/2},
\een
given by
\ben
(\tss x_{ij}\ |\ 1\leqslant i\leqslant m,\quad 1\leqslant j\leqslant n)
\mapsto
(\tss \la^{(r)}_{il}\ |\ 1\leqslant i,l\leqslant m,\quad 0\leqslant r\leqslant p)
\een
with $i>l$ for $r=0$,
thus completing the proof.
\epf

The argument used in the proof of Theorem~\ref{thm:secoadet}
provides a proof of the Poincar\'e--Birkhof\/f--Witt theorem for the algebra
$\wt\Y_q(\gl_m)$ and hence for $\Y_q(\gl_m)$.

\bco\label{cor:pbw}
Given any ordering on the set of generators of the algebra $\Y_q(\gl_m)$,
the ordered monomials in the generators form a basis of $\Y_q(\gl_m)$.
\eco

\bpf
As we showed in the proof of Theorem~\ref{thm:secoadet},
given any f\/inite family of ordered monomials in the generators
$\tau_{ij}^{(r)}$ of the algebra $\wt\Y_q(\gl_m)$,
we can choose a suf\/f\/iciently large value of $n$ such that
the images of the monomials under the homomorphism \eqref{wtyu}
are linearly independent. As the ordered monomials
span $\wt\Y_q(\gl_m)$, we may conclude that they form a basis of
$\wt\Y_q(\gl_m)$. The corresponding statement for the algebra $\Y_q(\gl_m)$
is now immediate from the relations~\eqref{taudef}
and~\eqref{reztze}.
\epf

We believe an analogue of the tensor product decomposition
\eqref{tensor} for the algebra $\Ar_m$ takes place. Namely,
denote by $\wt{\Ar}_0$ the commutative
subalgebra of $\Ar_m$ generated by the coef\/f\/icients of the virtual quantum
determinant $\qdet'(\Tc-\overline \Tc\tss u^{-1})^{}_{\Bc\Bc}$
with $\Bc=\{m+1,m+2,\dots\}$.

\bcj\label{conj:tenpr}
We have the tensor product decomposition
\beql{tproaymod}
\Ar_m=\wt{\Ar}_0\ot\wt\Y_q(\gl_m),
\eeq
where the Yangian $\wt\Y_q(\gl_m)$ is identified with its image under
the embedding $\vp$.
\ecj

\bre\label{rem:gln}
A slight modif\/ication of the proof of Theorem~\ref{thm:secoadet}
yields an alternative method to prove that the Yangian $\Y(\gl_m)$
is embedded into the Olshanski algebra $\Ar_m$; cf.~\cite{m:yt,o:ri}.
Moreover, as with the $q$-Yangian, this also provides
a new proof of the Poincar\'e--Birkhof\/f--Witt theorem for the Yangian
$\Y(\gl_m)$; cf. \cite{bk:pp,mno:yc,n:yq}.
\ere

\subsection*{Acknowledgements}

We would like to thank Boris Feigin for
discussions and Serge Ovsienko for help
with the proof of Theorem~\ref{thm:secoadet}.

\LastPageEnding

\end{document}